
\magnification \magstep1
\input amstex
\documentstyle{amsppt}
\rightheadtext{Uniformly Levi degenerate CR manifolds} \topmatter
\leftheadtext{Peter Ebenfelt}


\define \im{\text{\rm Im }}
\define \re{\text{\rm Re }}
\define \tr{\text{\rm tr }}
\define \ad{\text{\rm ad }}

\define \bR{\Bbb R}

\define \bC{\Bbb C}



\def\dim {\text {\rm dim}}

\def\Aut{\text{\rm Aut}}
\def\pb11{\overline{p_{11}}}
\def\pb02{\overline{p_{02}}}

\title Uniformly Levi degenerate CR manifolds; the 5
dimensional case
\endtitle
\author Peter Ebenfelt\footnote{Supported in part by a grant from the
Swedish Natural Science Research
Council.\newline}
\endauthor
\address Department of Mathematics, Royal Institute of Technology, 100
44 Stockholm, Sweden\endaddress
\email ebenfelt\@math.kth.se\endemail

\abstract In this paper, we consider real hypersurfaces $M$ in $\bC^3$
(or more generally, 5-dimensional CR manifolds of hypersurface type)
at   
uniformly Levi degenerate points, i.e.\ Levi degenerate points such
that the rank of the Levi form is constant in a neighborhood. We also
require the hypersurface to satisfy a certain second order
nondegeneracy condition (called $2$-nondegeneracy) at the point. For a 
real-analytic everywhere Levi degenerate hypersurface $M$ in $\bC^3$
which is not locally equivalent to a hypersurface of the form $\tilde
M\times \bC$, such points are dense
on $M$.

Our first result is the construction, near any point $p_0\in M$
satisfying the above conditions, of a principal bundle
$P\to M$ and a $\bR^{\dim P}$-valued $1$-form $\underline{\omega}$, 
uniquely determined by the CR structure on $M$, which defines an absolute
parallelism on $P$ and has the following property: Let $M$ and $M'$ be
two real-analytic hypersurfaces in $\bC^3$ with distinguished points
$p_0\in M$, $p_0'\in M'$ and parallelized principal bundles
$P,\underline{\omega}$, $P',\underline{\omega}'$, respectively. Then
there exists a local biholomorphism $h\:(\bC^3,p_0)\to (\bC^3,p_0')$
with $H(M)\subset M'$ if and only if there exists a real-analytic
diffeomorphism $H\:P\to P'$ with
$H^*\underline{\omega}'=\underline{\omega}$. ($H$ 
is then the lift of $h$, i.e. $\pi'\circ H=\pi\circ h$ where $\pi$ and
$\pi'$ denote the projections $\pi\:P\to M$, $\pi'\:P'\to M'$). This
solves the biholomorphic equivalence problem for uniformly Levi
degenerate hypersurfaces in $\bC^3$ at $2$-nondegenerate points in
view of Cartan's solution of the equivalence problem for absolute
parallelisms.  

A basic example of a hypersurface of the type under consideration is
the tube $\Gamma_\bC$ over the light cone. Our second result is the
characterization of 
$\Gamma_\bC$ by vanishing curvature conditions in the spirit of the
characterization of the unit sphere as the flat
model for strongly pseudoconvex hypersurfaces in $\bC^{n+1}$ in terms
of the Cartan--Chern--Moser connection.

\endabstract


\endtopmatter
\document

\heading 0. Introduction\endheading

\subhead 0.1. A brief history\endsubhead A fundamental problem in the
study of real submanifolds in complex 
space is the {\it biholomorphic equivalence problem} which in its most
general form asks for (intrinsic) conditions on two submanifolds
$M,M'\subset \bC^N$ at distinguished points $p_0\in M$, $p_0'\in M'$
which guarantee that there exists a local biholomorphism
$H\:\bC^N\to \bC^N$ defined near $p_0$ such that $H(p_0)=p_0'$ and
$H(M\cap U)= M'\cap U'$, for some open neighborhoods
$U,U'\subset \bC^N$ of $p_0$ and $p_0'$ respectively. When $M$ and
$M'$ are real-analytic, an equivalent formulation is to ask for a
real-analytic local CR diffeomorphism $f\:M\to M'$ defined near $p_0\in M$
with $f(p_0)=p_0'$. (For standard definitions and results on real
submanifolds in complex space and abstract CR structures, the reader
is referred e.g.\ to [BER]).       

The case where
$M$ and $M'$ are real-analytic and Levi
nondegenerate hypersurfaces was solved by Cartan [C1--2] in $\bC^2$,
and by Tanaka [T1--2] and Chern--Moser [CM] in $\bC^N$, $N\geq 2$. The
solution consists of producing a fiber bundle $Y\to M$, for
any given Levi nondegenerate hypersurface $M\subset \bC^N$, and 1-form
$\underline{\omega}$ on $Y$, valued in $\bR^{\dim Y}$, which at every
$y\in Y$ gives an isomorphism between $T_y Y$ and $\bR^{\dim Y}$ (an
absolute parallelism or $\{1\}$-structure on $Y$; see e.g. [KN] or [K]) such
that the following holds. If there
exists a CR diffeomorphism $f\:M\to M'$, then 
there exists a diffeomorphism $F\:Y\to Y'$ (where corresponding
objects for $M'$ are denoted with ') such that
$F^*\underline{\omega'}=\underline{\omega}$ and the following
diagram commutes
$$
\CD 
Y @> F>> Y'\\ @V \pi VV @VV \pi' V\\ M @>> f> M'.
\endCD\tag 0.1.1
$$ 
Conversely, if there exists a diffeomorphism $F\:Y\to Y'$ such that
$F^*\underline{\omega'}=\underline{\omega}$, then there exists a CR
diffeomorphism $f\:M\to M'$ such that \thetag{0.1.1} commutes. Suppose
that such a bundle $Y\to M$ and $\bR^{\dim Y}$-valued $1$-form
$\underline{\omega}$ can be constructed 
for every $M$ in some given class of manifolds. Then we shall
say that the bundle $Y\to M$ with $1$-form
$\underline{\omega}$ 
{\it reduces the CR structure on $M$ to a 
parallelism} (in this class). The
construction of a bundle 
$Y\to M$ which reduces the CR structure on $M$ to a parallelism for a
class of CR submanifolds in $\bC^N$ reduces the biholomorphic
equivalence problem for the real-analytic manifolds in this class to
the equivalence problem for $\{1\}$-structures. The latter problem was
solved by 
Cartan, and is well understood (see e.g. [G] or [S]). 

The bundle $Y\to M$ constructed in [CM] is in fact a principal fiber
bundle with group 
$G_0$, where $G_0$ is the isotropy subgroup of $SU(p+1,q+1)$, 
$p+q=N-1$, and $p$, 
$q$ are the number of positive and negative eigenvalues,
respectively, of the Levi form. The authors of [CM] also construct a Cartan
connection $\Pi$ valued in the Lie algebra $\frak s\frak u(p+1,q+1)$ which
defines the same parallelism as $\underline{\omega}$ (given a suitable
identification of $\bR^{\dim Y}$ with $\frak s\frak u(p+1,q+1)$). Covariant 
differentiation of the curvature $\Omega:=d\Pi-\Pi\wedge\Pi$ produces a
complete set of invariants for a real-analytic Levi nondegenerate
hypersurface. In particular, it follows that a real-analytic strongly
pseudoconvex hypersurface in $\bC^N$ is locally biholomorphic to a
piece of the $(2N-1)$-sphere (the ``standard model'' for such
hypersurfaces) if and only if the curvature $\Omega$ is identically
zero (i.e.\ the connection is flat). The reader is also referred to
the work of Burns-Shnider [BS] and Webster [We] for further
discussion in the Levi nondegenerate case. 

More recently, CR manifolds
of higher codimension whose Levi forms are suitably nondegenerate have
been studied by several authors. Since the main focus in the present
paper is on hypersurfaces, we mention only the papers by \v
Cap--Schichl [CS], Ezhov--Isaev--Schmalz [EIS], Garrity--Mizner [GM],
Schmalz--Slov\'ak 
[SS], and refer the interested reader to these papers for further
information about the higher codimensional case. 

In this paper, we consider
real hypersurfaces (and, more generally, CR 
manifolds of hypersurface type) which have {\it degenerate} Levi
forms. Before describing our main results, we
should mention that another approach to the biholomorphic equivalence
problem is via normal forms. Normal forms for certain types of
Levi degeneracies were studied by the author in [E3--4]; another class
of Levi degeneracies in $\bC^2$ was considered by Wong [Wo]. However,
at least to the best of the author's knowledge, the geometric approach
as described above has not been previously studied for CR manifolds
with degenerate Levi forms. The idea in this paper is to use the
higher order invariant tensors introduced by the author in [E4] as a
complement to the degenerate Levi form. We mostly restrict our
attention to 5-dimensional 
manifolds in order to keep the number of cases and the notation to a
minimum. The main results can be generalized to higher dimensional
manifolds under Conditions 2.21 and 2.25; see the concluding remarks
in \S 5. 

The paper is organized as follows. Our main results are explained in
\S0.2. In \S0.3, some examples of everywhere Levi degenerate hypersurfaces
that arise e.g.\ in PDE theory are given. \S1 is devoted to
preliminary material including basic definitions and properties of
Levi uniform CR manifolds. The necessary constructions for the main
results are given in \S2--3, and in \S 4 a discussion and
characterization of the tube over the light cone is given. The paper
concludes with some remarks in \S 5 about the higher dimensional case. 

\subhead 0.2. The main results\endsubhead Our main results concern
real hypersurfaces $M$ in 
$\bC^3$ or, more generally, 5-dimensional CR manifolds of hypersurface
type, which are uniformly Levi degenerate in 
the sense that the Levi form has one nonzero and one zero eigenvalue
in a neighborhood of a distinguished point $p_0\in M$ ({\it Levi
uniform} of rank $1$ according to Definition 1.6 below). We also require
$M$ to be $2$-nondegenerate at $p_0$ (see section 1 or [BER, Chapter
XI]). The latter condition guarantees that if $M$ is a real hypersurface
in $\bC^3$, then it is {\it holomorphically
nondegenerate} (see [BER, Chapter XI]) and, in particular, not locally
biholomorphic to 
a manifold of the form $\tilde M\times \bC$ where $\tilde M$ is a
hypersurface in $\bC^2$. (However, $M$ is always foliated by complex
lines, but this foliation cannot be ``straightened''; see Proposition
1.15.) In fact, if $M$ is real-analytic and everywhere Levi
degenerate, then at most points $p$ (off a proper real-analytic
subvariety) $M$ is either locally biholomorphic to $\tilde M\times \bC$ for
some $\tilde M\subset \bC^2$ or $M$ is Levi uniform of rank $1$ and
$2$-nondegenerate at $p$.  

The most  
important example (indeed, the 
``standard model'' for such manifolds) is the tube in
$\bC^3$ over 
the light cone in $\bR^3$ (see Example 1.7 and section
4 for further discussion); other examples of everywhere Levi
degenerate hypersurfaces that arise naturally are  
given, for motivation, in the next section. Let us just point out that
the biholomorphically invariant geometry of the tube $\Gamma_\bC$  over 
the light cone in $\bC^4$ plays an important role in e.g.\
axiomatic quantum 
field theory since $\Gamma_\bC$ bounds the so-called past and future
tubes; see e.g.\ Sergeev--Vladimirov [SV] and Zhou [Z] and the references in
these papers.

For the class of
hypersurfaces described above, we define a new CR invariant $\hat k$ (see
\thetag{2.34}); more generally, for
$(2n+1)$-dimensional Levi uniform CR hypersurfaces of rank $n-1$ which
satisfy Conditions 2.21 and 2.25, we introduce a sequence of 
invariants which in the 5-dimensional case reduce to the single
invariant $\hat k$. One of our main results is the following. We refer
the reader to section 1 for relevant definitions.

\proclaim{Theorem 1} Let $M$ be a $5$-dimensional CR manifold of
hypersurface type which is $2$-nondegenerate and Levi uniform of rank
$1$ at $p_0\in M$. Then, there exists a principal fiber bundle $P\to
M$ with group $G_0$ and a $1$-form $\underline{\omega}$ on $P$ which
defines an isomorphism between $T_u P$ and $\bR^{\dim P}$ for every
$u\in P$ and reduces the CR structure on $M$ to a
parallelism. The group $G_0$ is a subgroup of $GL(\bR^{3})$ which has
dimension two if the invariant $|\hat k(p_0)|=2$ and dimension one
otherwise. \endproclaim 

Theorem 1 will be a consequence of the more detailed Theorems
3.1.37 and 3.2.9. We should mention that the group $G_0$, and hence
the bundle 
$P$, in Theorem 1 is disconnected and has two components. In order to
obtain a connected bundle, we have to choose an ``orientation'' for the
Levi nullspace as explained \S$2$ (see Theorems 3.1.37 and 3.2.9). 

As mentioned above, the most important example is the tube
$\Gamma_\bC$ over the light cone for which the invariant $\hat k\equiv
2i$. We shall now characterize $\Gamma_\bC$ among all $M$, as in Theorem
1 with $\hat k\equiv 2i$, by a curvature condition in the spirit of the
characterization of the sphere among strongly pseudoconvex
hypersurfaces as described section 0.1. 
There is a subgroup $H$ of $GL(\bC^4)$ and a subgroup $H_0$ of $H$
such that $H$ can be viewed as a principal fiber bundle over
$\Gamma_\bC$ with group $H_0$. The matrix valued Maurer-Cartan forms
$\Pi$ of $H$ 
define a Cartan connection on $\Gamma_\bC$ valued in $\frak h$, the
Lie algebra of $H$, with
vanishing curvature $\Omega=d\Pi-\Pi\wedge\Pi$. (All this is explained
in detail in section 4.)
For a real-analytic CR manifold $M$ as in Theorem 1 with
the invariant $\hat k\equiv 2i$, we can identify the group $G_0$ with
the group $H_0$, and construct, using $\underline{\omega}$, a $\frak
h$-valued 1-form $\Pi$ which, unfortunately, in general is not a
Cartan connection. However, we have the following result, which is a
consequence of the more detailed Theorem 4.31. 

\proclaim{Theorem 2} Let $M$ a real-analytic CR manifold satisfying
the conditions in Theorem $1$ with $\hat k\equiv
2i$. Then there exists a $\frak h$-valued $1$-form $\Pi$ on the
principal bundle $P\to M$, given by Theorem $1$, which gives an
isomorphism between $T_u P$ and $\frak h$ for every $u\in P$ and with
the following property. There exists a real-analytic CR diffeomorphism
$f\: M\to \Gamma_\bC$, defined near $p_0\in M$, if and only if the
curvature 
$$
\Omega:=d\Pi-\Pi\wedge\Pi\tag 0.2.1
$$ 
vanishes identically near $p_0$.\endproclaim

We conclude this section by giving an application of Theorem 1. Let
$\Aut(M,p_0)$ denote the stability group of a CR manifold $M$ at
$p_0\in M$, i.e. the group of germs at $p_0$ of local smooth CR
diffeomorphisms $f\:(M,p_0)\to (M,p_0)$. 
Suppose that $M$ satisfies
the conditions of Theorem 1 at $p_0$. Pick any point $u\in P_{p_0}$,
where $P\to M$ is the principal $G_0$ bundle given by Theorem 1 and
$P_p$ denotes the fiber over $p\in M$. By [K, Theorem 3.2] and
Theorem 1, the group $\Aut(M,p_0)$ embeds as a closed submanifold of
the fiber $P_{p_0}\cong G_0$ via the mapping
$$
\Aut (M,p_0)\ni f\mapsto F(u)\in P_{p_0},\tag 0.2.2
$$
where $F\:P\to P$ is the lift of $f$ as in the diagram
\thetag{0.1.1}. Thus, $\dim \Aut(M,p_0)$ is at most $2$ if the
invariant $|\hat k(p_0)|=2$ and at most $1$ if $|\hat k(p_0)|\neq
2$. We formulate this as follows.

\proclaim{Corollary 3} Let $M$ be a $5$-dimensional CR manifold of
hypersurface type which is $2$-nondegenerate and Levi uniform of rank
$1$ at $p_0\in M$. Then $\dim \Aut(M,p_0)\leq 2$.
\endproclaim

The bound in Corollary 3 cannot be improved, since
$\dim\Aut(\Gamma_\bC,p)=2$ for any $p\in \Gamma_\bC$ as is shown in
section 4. We should mention that in the recent preprint [Er], it was shown
that the bound $\Aut (M,p_0)\leq 3$ holds for the class of all
real-analytic $2$-nondegenerate hypersurfaces in $\bC^3$; observe that
when $M$ is real-analytic then, by the reflection 
principle (see [BJT]), every $f\in \Aut(M,p_0)$ is real-analytic,
since $2$-nondegeneracy implies essential finiteness (see [BER,
Chapter XI]). 

\subhead 0.3. Examples of everywhere Levi degenerate CR
manifolds\endsubhead Trivial examples of real
hypersurfaces in $\bC^{n+1}$ which are everywhere Levi degenerate can
be obtained by taking any hypersurface of the form $\tilde M\times
\bC$, where $\tilde M$ is a real hypersurface in $\bC^n$. Such
hypersurfaces, as mentioned in the previous section, are never
$2$-nondegenerate and, hence, are not of interest to us in the present
paper. We give here
two (from our viewpoint) more interesting classes of everywhere Levi
degenerate hypersurfaces in $\bC^{n+1}$. The reader is referred to
section 1 for relevant definitions. 

\example{Example 0.3.1 (Everywhere characteristic hypersurfaces)} Let
$p(x)$ be a homogeneous polynomial of $x=(x_1,\ldots, x_{n+1})$. A
real hypersurface $M\subset\bC^{n+1}$ is called {\it characteristic}
at $p_0\in M$ for the partial differential operator
$p(\partial):=p(\partial/\partial Z)$ if 
$$
p\left(\partial \rho(p_0,\bar p_0)/\partial Z\right )=0,\tag
0.3.2
$$
where $\rho(Z,\bar Z)=0$ is a defining equation for $M$ near
$p_0$ and $$\partial/\partial Z=(\partial/\partial
Z_1,\ldots,\partial/\partial Z_{n+1}).\tag 0.3.3
$$
$M$ is called {\it everywhere
characteristic} if $M$ is 
characteristic at every point. Everywhere characteristic hypersurfaces
for a given operator $p(\partial)$ arise as natural boundaries for the
holomorphic continuation of (holomorphic) solutions of 
$p(\partial)
u=0$ (see e.g.\ [H\"o, Chapter IX.4]).
A concrete example is given by the so-called {\it Lie ball} defined by
the equation
$$
|Z|^2+\left(|Z|^4-\left|\sum_{k=1}^{n+1}Z_k^2\right|^2\right)^{1/2}<1.\tag
0.3.4
$$
The Lie ball is the maximal domain in $\bC^{n+1}$ to which every
harmonic function in the 
unit ball of $\bR^{n+1}$ can be holomorphically continued (see
e.g. [A]; cf.\ also [E1]). The boundary
of the Lie ball
is everywhere characteristic (at every smooth point) for the ``Laplace
operator'' 
$$
\sum_{j=1}^{n+1}(\partial/\partial Z_k)^2.\tag 0.3.5
$$ 
Another example is
the tube over the light cone (see Example 1.7) which is everywhere
characteristic for the ``wave operator'' 
$$
\sum_{j=1}^{n}(\partial/\partial
Z_k)^2-(\partial/\partial Z_{n+1})^2.\tag 0.3.6
$$
We have the following.
\proclaim{Proposition 0.3.7} Let $M\subset\bC^{n+1}$ be a real
hypersurface which is everywhere
characteristic for a homogeneous partial differential operator
$p(\partial)$. Then $M$ is everywhere Levi degenerate. 
\endproclaim

\demo{Proof} Pick $p_0\in M$ and let $\rho(Z,\bar Z)=0$ be a
defining equation for $M$ near $p_0\in M$. 
We first claim that the CR vector field 
$$
L:=\sum_{k=1}^{n+1} \overline{\partial p(\partial
\rho/\partial Z)/\partial
x_k}\frac{\partial}{\partial \bar Z_k} 
\tag 0.3.8
$$
is tangent to $M$. Indeed, since $M$ is everywhere characteristic for 
$p(\partial)$, we have
$$
p\left(\partial\rho/\partial Z\right)=a\rho\tag 0.3.9
$$ for some function $a$. The claim now follows from Euler's formula. 
By differentiating \thetag{0.3.9}, it is straightforward (and left to
the reader) to verify
that $L$ is a nullvector for the 
Levi form at every $p\in M$ near $p_0$. This proves the
proposition.\qed \enddemo
\endexample

\example{Example 0.3.10} Let $p(x)$ be a homogeneous polynomial in
$x=(x_1,\ldots, x_{n+1})$ with real coefficients and assume that
$\partial p/\partial x$ is not identically zero along the variety
$V_\bR:=\{x\in\bR^{n+1}\: p(x)=0\}$. Then the tube $V_\bC$ over $V_\bR$
in $\bC^{n+1}$, 
$$
V_\bC:=\{Z\in \bC^{n+1}\:p(\re Z)=0\},\tag 0.3.11
$$ 
is a real hypersurface (outside a lower dimensional real algebraic
subvariety) which we denote by $M$. The ``radial''
CR vector field 
$$
L=\sum_{j=1}^{n+1}\re Z_j\, \partial/\partial \bar Z_j\tag 0.3.12
$$ is tangent to $M$, and the reader can easily verify that $L$ is a
null vector for the Levi form of $M$ at every $p\in M$. A concrete
example is again the tube over the light cone (Example 1.7). Another
example is the cubic defined by \thetag{0.3.11} with
$p(x)=x_1^3+x_2^3-x_3^3$ which was given by Freeman [F] as an example
of a manifold foliated by complex curves but not locally biholomorphic
to a manifold of the form $\tilde M\times \bC$.
\endexample

\heading 1. Preliminaries\endheading

Let $M$ be a CR manifold with  CR bundle  $\Cal V$. Recall that
this means that $\Cal V$ is a subbundle of the compexified tangent
bundle $\bC TM$ such that $\Cal V_p\cap \bar\Cal V_p=\{0\}$ for
every $p\in M$, and $\Cal V$ is {\it formally integrable} i.e. any
commutator between sections of $\Cal V$ is again a section of
$\Cal V$; sections of $\Cal V$ will henceforth be called {\it CR
vector fields}. We shall denote the CR dimension of $M$, i.e. the
(complex) dimension of the fibers $\Cal V_p$ for $p\in M$, by $n$.
We shall assume that $M$ is of {\it hypersurface type},
i.e. the complex
dimension of $T^0_pM:=(\Cal V_p\oplus\bar \Cal
V_p)^\perp\subset\bC T^*_p M$, for $p\in M$, is one. In
particular, the dimension of $M$ is $2n+1$. For the remainder of this
paper unless explicitly stated otherwise, all CR manifolds will be of
hypersurface type.
The bundle $T^0
M\subset\bC T^* M$ is called the {\it characteristic bundle}, and
real sections of $T^0M$ are called {\it characteristic forms}. The
subbundle $T'M\subset\bC T^*M$, defined at $p\in M$ by $T'_p
M:=\Cal V_p^\perp$, is called the {\it holomorphic cotangent
bundle}. The formal integrability of $\Cal V$ is
equivalent to the following property of $T'M$: If $\omega$ is a
section of $T'M$, then $d\omega$ is a section of the ideal generated
by $T'M$ in the exterior algebra of $\bC TM$.

Let $L_{\bar 1},\ldots, L_{\bar n}$ be a basis for the CR vector
fields near some distinguished point $p_0\in M$. Also, let
$\theta$ be a nonvanishing characteristic form near the same point
$p_0$. Following [E2] (see also [BER] and [E4]), we define linear
operator $\Cal T_{\bar 1},\ldots \Cal T_{\bar n}$ on the holomorphic
1-forms on
$M$, i.e. the sections of $T'M$, near $p_0$ as follows
$$
T_{\bar A}\omega:= \frac{1}{2i} L_{\bar A}\lrcorner d\omega,\tag
1.1
$$
where $\lrcorner$ denotes the usual contraction by a
vector field. For $p\in M$ near $p_0$ and positive integers $k$,
we define the subspace $E_{k,p}\subset T'_pM$ as the (complex)
linear span of $\theta_p$ and $(T_{\bar A_j}\ldots
T_{\bar A_1}\theta)_p$, for all $1\leq j\leq k$ and all $j$-tuples
$(A_1,\ldots,A_j)\in\{1,\ldots, n\}^j$. We define
$E_{0,p}$ to be $T^0_pM$. The CR manifold $M$ is said to be {\it
finitely nondegenerate} at $p$ of
$$
E_{k,0}=T'_p M\tag 1.2
$$
for
some integer $k\geq 1$, and $k_0$-nondegenerate at $p$ if $k_0$ is
the smallest integer $k$ for which \thetag{1.2} holds. It was
shown in [E2] (see also [BER, Chapter XI]) that this definition is
consistent with the one for real hypersurfaces of
$\bC^{n+1}$ given in [BHR]. (These notions can also be extended to
CR manifolds of arbitrary codimension; see [BER] or [E4].) For
each integer $k$ such that $E_k(p)\neq T'_pM$, the  author
introduced in [E4] an invariant tensor
$$
\psi_{k+1}\in \Cal
V_p^*\otimes\ldots\Cal V_p^*\otimes F_{k,p}^*\otimes (T^0 _pM)^*,\tag
1.3
$$
where $\Cal V_p^*$ occurs
$k$ times in \thetag{1.3} and $F_{k,p}=E_{k,p}^\perp\cap \bar \Cal
V_p$. The sequence of tensors $\psi_2,\ldots, 
\psi_{k_0+1}$ describes
in more detail the data associated with $k_0$-nondegeneracy. In
this paper, we shall mainly consider certain classes of
2-nondegenerate CR manifolds and, hence, we shall only be interested
in the first two tensors $\psi_2$ and $\psi_3$. The reader is referred
to [E4] for the precise definition of the tensors $\psi_k$ and their
basic properties.

The second order tensor $\psi_2$ at a point $p\in M$ is given in
the bases $L_{\bar 1,p},\ldots, L_{\bar n,p}$ of $\Cal V_p$ and
$\theta_p$ of $T^0_p M$ as an $n\times n$ matrix $(g_{\bar A
B}(p))_{1\leq A,B,\leq n}$, where $$ g_{\bar A
B}(p):=\left<(T_{\bar A}\theta)_p,L_{B,p}\right>.\tag 1.4 $$ We
use here, and throughout this paper, the convention that
$L_B=\overline{L_{\bar B}}$. The reader should observe that we
have the identity $$ g_{\bar A
B}(p)=\frac{1}{2i}\left<\theta_p,[L_B, L_{\bar A}]_p\right>, \tag
1.5 $$ which coincides with the {\it Levi form} $\Cal L_\theta$
of $M$ at $p$ and $\theta_p$. From this observation, we see that
$1$-nondegeneracy of $M$ at $p$ is equivalent to the classical
notion of {\it Levi nondegeneracy} of $M$ at $p$. Moreover, the
subspace $F_{1,p}\subset \bar\Cal V_p$ coincides with the {\it
Levi nullspace} of the Levi form, i.e. those vectors $X_p\in
\bar\Cal V_p$ for which the linear form $Y_p\mapsto \Cal
L_\theta(Y_p,X_p)$ is zero. For the remainder of this paper, we
shall use the notation $\frak N_p$ for the Levi nullspace
$F_{1,p}$.

\definition{Definition 1.6} A CR manifold $M$ of hypersurface type is
{\it $1$-uniform} or {\it Levi uniform} (of rank $r$) at $p\in M$
if the rank of the Levi form is constant (and equal to $r$) in a
neighborhood of $p$.
\enddefinition

Observe that the extremal cases of CR manifolds which are Levi
uniform of rank 0 or $n$ at a point $p$ are precisely those which
are Levi flat or Levi nondegenerate, respectively, at $p$. Such CR
manifolds which, in addition, are real-analytic are by now fairly
well understood: a real-analytic Levi flat CR manifold is locally
CR equivalent to the real hyperplane $\im Z_{n+1}=0$ in
$\bC^{n+1}$, and a theory for real-analytic Levi nondegenerate CR
manifolds was developed by E. Cartan [C1-2], Tanaka [T1-2],  and
Chern--Moser [CM].

The reader should also observe that any real-analytic CR manifold
is Levi uniform outside a proper real-analytic subvariety (in
particular, on a dense open subset). Before proceeding, let us
pause and give an example of a Levi uniform CR manifold which is neither Levi
nondegenerate nor Levi flat. 

\example{Example 1.7} The tube in $\bC^3$ over the light cone in
$\bR^3$, i.e. the variety defined by $$ (\re Z_1)^2+(\re
Z_2)^2-(\re Z^3)^2=0,\tag 1.8 $$ is Levi uniform of rank 1 at every
nonsingular point, i.e. at every point where it is a real
submanifold. Thus, it is Levi
uniform, but neither Levi flat nor Levi nondegenerate since the CR
dimension $n$ is $2$. The reader can also
verify that the CR manifold given by equation \thetag{1.8} is
$2$-nondegenerate at every nonsingular point. This example will be
discussed in greater detail in \S4 below. (See also [E3],
where this example is further discussed in connection with a normal
form for $2$-nondegenerate hypersurfaces in $\bC^3$.)
\endexample

Note that if $M$ is Levi uniform (of rank $r$) at $p_0\in M$, then
 the
subspaces $\frak N_p$ for $p$ near $p_0$ form a (rank $n-r$)
subbundle $\frak N$ of $\bar \Cal V$. From now on, we assume that
$M$ is Levi uniform of rank $r$, with $0<r<n$,  in a neighborhood
of $p_0$ (to which we restrict our attention). We may arrange our
basis for the CR vector fields $L_{\bar 1}, \ldots, L_{\bar n}$ so
that $L_{r+1},\ldots,L_n$ is a basis for the sections of $\frak N$
near that point. The second order tensor for $p$ near $p_0$ then
takes the form $$ (g_{\bar A B})=\pmatrix g_{\bar \alpha\beta} &
0\\ 0& 0\endpmatrix, \tag 1.9 $$ where $(g_{\bar
\alpha\beta})_{1\leq \alpha,\beta\leq r}$ is an $r\times r$
nondegenerate Hermitian matrix of smooth functions. In what
follows, we shall use the summation convention and also the
convention that capital roman indices $A,B,\ldots$ run over the
integers $\{1,\ldots,n\}$ and Greek indices $\alpha,\beta,\ldots$
run over $\{1,\ldots,r\}$. The third order tensor $\psi_3$ can be
represented near $p_0$ by $n\times n$ matrices $(h_{\bar A\bar B
k})_{1\leq A,B\leq n}$, $k=r+1,\ldots, n$, of smooth functions,
where $$ h_{\bar A\bar B k}:=\left<\Cal T_{\bar B}\Cal T_{\bar A}
\theta,L_{k}\right>.\tag 1.10 $$

\proclaim{Proposition 1.11} Assume that $M$ is Levi uniform of rank
$r$ at
$p_0$. Then, in the notation introduced
above, for every $k=r+1,\ldots, n$ and  all
$p$ in a neighborhood of $p_0$, it holds that $h_{\bar A\bar B
k}=0$ whenever $A$ or $B$ belongs to
$\{r+1,\ldots, n\}$.
\endproclaim

\demo{Proof} Using a well known identity (see e.g.\ [He, Chapter
I.2]; see also the remark concerning our normalization of the
pairing $\left<\cdot,\cdot\right>$ in [E4]), we have $$\aligned
h_{\bar A\bar B k} :&= \left<\Cal T_{\bar B}\Cal T_{\bar
A}\theta,L_k\right>=\left< d(\Cal T_{\bar A}\theta), L_{\bar
B}\wedge L_k\right>
\\& = L_{\bar B}(\left<\Cal T_{\bar A}\theta,L_k\right>)
-L_k(\left<\Cal T_{\bar A}\theta,L_{\bar B}\right>)-\left<\Cal
T_{\bar A}\theta,[L_{\bar B},L_k]\right>\\ & =-\left<\Cal T_{\bar
A}\theta,[L_{\bar B},L_k]\right>,
\endaligned\tag 1.12
$$ where the last identity follows from the facts that $\left<\Cal
T_{\bar A}\theta,L_k\right>\equiv 0$ and $\left<\Cal T_{\bar
A}\theta,L_{\bar B}\right>\equiv 0$. It is proved in [E4] that the
matrices $h_{\bar A\bar B k}$ are symmetric, so to prove the
proposition it suffices to show $h_{\bar A\bar l k}=0$ in a
neighborhood of $p_0$ for $k,l\geq r+1$, i.e.\ $\left<\Cal T_{\bar
A}\theta,[L_{\bar l}, L_k]\right>=0$ in view of \thetag{1.12}. To
this end note, using the fact that $L_{k,p}$ and $L_{l,p}$ are
null vectors for the Levi form at every $p$ near $p_0$, that $$
\left<\Cal T_{\bar A}\theta,[L_{\bar l},
L_k]\right>=-\left<\theta,[L_{\bar A},[L_{\bar l},
L_k]]\right>.\tag 1.13 $$ Thus, by also using the Jacobi
identity, we obtain $$\left<\Cal T_{\bar A}\theta,[L_{\bar l},
L_k]\right>=\left<\theta,[L_{\bar l},[ L_k,L_{\bar A}]]\right>+
\left<\theta,[L_k,[L_{\bar A},L_{\bar l}]]\right>.\tag 1.14$$ The
second term on the right hand side of \thetag{1.14} vanishes
since $[L_{\bar A},L_{\bar l}]$ is a CR vector field by the formal
integrability of $\Cal V$ and $L_k$ is a null vector field for the
Levi form. To show that the first term also vanishes, we must show
that $[ L_k,L_{\bar A}]$ is a section of $\Cal V\oplus\bar \Cal
V$. This fact follows again from the fact that $L_k$ is a null
vector field for the Levi form since the latter is equivalent to
$\left<\theta,[L_{ k},L_{\bar A}]\right>=0$ for every $A=1,\ldots,
n$. The proof of Proposition 1.11 is complete. \qed
\enddemo

Let us digress
briefly to note
the following result which, although of no importance for the
remainder of this paper, follows from (the proof of)
Proposition 1.11 above.

\proclaim{Proposition 1.15} If $M\subset \bC^{n+1}$ is a real
hypersurface which
is Levi uniform of rank $r<n$ at $p_0\in M$, then $M$ is
foliated by complex manifolds of dimension $n-r$ in a
neighborhood of $p_0$. \endproclaim

\remark {Remarks $1.16$} \roster
\item"(i)" The following partial converse to Proposition
1.15 is easy to verify: if $M$ is foliated by complex manifolds of
dimension $n-r$ then the rank of the Levi form at any point is $\leq
r$.
\item"(ii)" We should point out that the foliation given by
Proposition 1.15, even when $M$ is real-analytic, can in general not
be ``straightened'', i.e. it is
in general not true that $M$, as in the proposition, is CR equivalent
to a real hypersurface of the form $\tilde M\times \bC^r\subset
\bC^{n+1}$, where $\tilde M$ is a real hypersurface in
$\bC^{n+1-r}$. Indeed, if $M$ is CR equivalent to $\tilde M\times
\bC^r$ (which is a {\it holomorphically degenerate} hypersurface; see
e.g. [BER, Chapter XI]), then it cannot be finitely nondegenerate at
any point. Thus, the foliation of the tube over the light cone
(Example 1.7) cannot be straightened.
\endroster
\endremark \medskip

\demo{Proof of Proposition $1.15$} An immediate consequence of
Proposition 1.11 and \thetag{1.12} is that, for $k,l\geq r+1$,
$$
[L_k,L_{\bar l}]=\sum_{m=r+1}^n (a^m_{k\bar l}L_m+b^{\bar m}_{k\bar
l}L_{\bar
m}),\tag 1.17
$$
where the $a_{k\bar l}^m$ and $b^{\bar m}_{k\bar l}$ are smooth
functions
satisfying $a^m_{l\bar k}+\overline{b^{\bar m}_{k\bar l}}=0$. Thus, by
the
Frobenius theorem, $M$ is
foliated near $p_0$ by $2(n-r)$-dimensional integral manifolds of $\re
L_k,\im L_k$, $k=r+1,\ldots, n$. By the Newlander--Nirenberg theorem
and
\thetag{1.17}, these manifolds are $(n-r)$-dimensional complex
submanifolds in $\bC^{n+1}$.\qed
\enddemo

Returning to the third order tensor $\psi_3$, we observe that
Proposition 1.11 shows that the matrices $(h_{\bar A\bar B k})$
representing $\psi_3$ are of the form
$$
(h_{\bar A\bar B k})=\pmatrix h_{\bar\alpha\bar\beta
k}&0\\0&0\endpmatrix.\tag 1.18
$$
We conclude this section with the following observation, whose proof
is immediate and left to the reader, characterizing $2$-nondegeneracy
for Levi uniform CR
manifolds of rank $r$ in terms of the $r\times r$ matrices
$(h_{\bar\alpha\bar \beta k})$ in \thetag{1.18}.

\proclaim{Proposition 1.19} Assume that $M$ is Levi uniform of rank
$r$ at $p_0$. Then, in the notation introduced above, $M$ is
$2$-nondegenerate at $p_0$ if and only if the symmetric matrices
$(h_{\bar\alpha\bar\beta k}(p_0))_{1\leq \alpha,\beta\leq r}$,
$k=r+1,\ldots, n$, are linearly independent over $\bC$.\endproclaim

\heading 2. A $G$-structure for Levi uniform
CR manifolds of rank $n-1$\endheading

We keep the notation from the previous section. However, from now
on, we shall restrict ourselves to the case $r=n-1$, i.e.\ the
case where the rank of the Levi form near $p_0$ is $n-1$. Thus, we assume
that $M$ is a smooth CR manifold (of hypersurface type) of CR
dimension $n$ which is Levi uniform of rank $n-1$ at a
distinguished point $p_0\in M$. We shall restrict our attention to
a small neighborhood of $p_0$. In what follows, $M$ will denote a
sufficiently small neighborhood of $p_0$. We have the following
invariant subbundles of the cotangent bundle $\bC T^* M$ $$ T^0
M\subset T'' M\subset T' M,\tag 2.1 $$ where $T^0 M$ and $T'M$
were introduced in \S 1 and where $T''M$ is defined by $$ T''_p
M=\{ \omega_p\in T'_p M\: \left < \omega_p, L_p\right>=0,\ \forall
L_p\in \frak N_p\}.\tag 2.2$$ Observe that $$ \dim T^0_p M=1,\
\dim T''_p M = n,\ \dim T'_p M = n+1.\tag 2.3$$ Let
$\theta,\theta^1,\ldots,\theta^n$ be a basis for the holomorphic
1-forms (i.e.\ sections of $T'M$) with the additional properties
that $\theta$ is real and a basis for the sections of $T^0 M$ and
$\theta,\theta^1,\ldots, \theta^{n-1}$ is a basis for the sections
of $T'' M$. Any other such basis $\tilde \theta,\tilde
\theta^1,\ldots,\tilde \theta^n$ is related to
$\theta,\theta^1,\ldots,\theta^n$ by $$ \pmatrix\tilde\theta\\
\tilde \theta^\alpha\\ \tilde \theta^n
\endpmatrix=\pmatrix u&0&0\\ u^\alpha & u^\alpha_\beta&0 \\ \xi^\alpha &
\xi^\alpha_\beta& \xi\endpmatrix \pmatrix \theta\\ \theta^\beta\\
\theta^n\endpmatrix,\tag 2.4 $$ where the coefficients in the
$(n+1)\times (n+1)$ matrix in \thetag{2.4} are smooth functions
(all complex valued, except $u$ which is real valued); also,
recall that we are using the summation convention, and Greek
indices run over the set $\{1,\ldots, n-1\}$ since $r=1$ here.

Observe that $T''_pM\cap\overline{T''_pM}=T^0_pM$ and
$T''M\cup\overline{T''M}$ is a rank $2n-1$ subbundle of $\bC T^*M$. The
$1$-forms $\theta,\theta^\alpha,\theta^{\bar\alpha}$, where
$\theta^{\bar\alpha}=\overline{\theta^\alpha}$ as mentioned in \S 1,
yield a coframe for bundle $T''M\cup\overline{T''M}$. Consider the
bundle  $Y\to M$ consisting of all such coframes
$(\omega,\omega^\alpha,\omega^{\bar\alpha})^\tau$,
$$ \pmatrix\omega\\ \omega^\alpha\\ \omega^{\bar \alpha}
\endpmatrix=\pmatrix u&0&0\\ u^\alpha & u^\alpha_\beta&0 \\
\overline{u^\alpha} &0& \overline{u^\alpha_\beta}
\endpmatrix \pmatrix \theta\\ \theta^\beta\\
 \theta^{\bar \beta}\endpmatrix,\tag
2.5 $$
where $(u,u^\alpha,u^\alpha_\beta)\in (\bR\setminus\{0\})\times\bC^{n-1}\times
GL(\bC^{n-1})$. If we let $G\subset GL(\bC^{2n-1})$ denote the
group consisting of matrices of the form $$S=\pmatrix u&0&0\\
u^\alpha & u^\alpha_\beta&0 \\ \overline{u^\alpha} &0&
\overline{u^\alpha_\beta}
\endpmatrix,\quad (u,u^\alpha,u^\alpha_\beta)\in (\bR\setminus\{0\})
\times\bC^{n-1}\times
GL(\bC^{n-1}),\tag 2.6$$
then $Y\to M$ is a principal fiber bundle over $M$
with group $G$; \thetag{2.5} gives a trivialization of $Y$ in which
$(u,u^\alpha, u^\alpha_\beta)$, or $S$ given by \thetag{2.6}, are
(global) coordinates of $Y$. 
We denote by $\frak g$ the Lie algebra of $G$,
i.e. the space of matrices
$$T=\pmatrix v&0&0\\
v^\alpha & v^\alpha_\beta&0 \\ \overline{v^\alpha} &0&
\overline{v^\alpha_\beta}
\endpmatrix,\quad (v,v^\alpha,v^\alpha_\beta)\in \bR
\times\bC^{n-1}\times
\Cal M(\bC^{n-1}),\tag 2.7$$
where $\Cal M(\bC^{n-1})$ denotes the space of all $(n-1)\times(n-1)$
matrices. If we pull back the forms $\theta,
\theta^A,\theta^{\bar A}$ (where capital Roman indices run over
$\{1,\ldots, n\}$) to $Y$, still denoting the
pulled-back forms by $\theta,\theta^A,\theta^{\bar A}$, then
\thetag{2.5} defines $1$-forms
$\omega,\omega^{\alpha},\omega^{\bar\alpha}$ on $Y$. The reader can
verify that the latter 1-forms are 
invariantly defined on $Y$, i.e. independent of the initial choice of
$\theta,\theta^{\alpha},\theta^{\bar\alpha}$ above.

Differentiating the 1-forms
$\omega,\omega^\alpha,\omega^{\bar \alpha}$ we obtain
$$
d\pmatrix\omega\\ \omega^\alpha\\ \omega^{\bar \alpha}
\endpmatrix=dS S^{-1}\wedge\pmatrix \omega\\ \omega^\beta\\
 \omega^{\bar \beta}\endpmatrix+Sd \pmatrix \theta\\ \theta^\beta\\
 \theta^{\bar \beta}\endpmatrix,\tag 2.8
$$
where $S\in G$ is given by \thetag{2.6}. The elements of the matrix
valued 1-form $dS S^{-1}$ are Maurer--Cartan forms for the Lie group
$G$ (see e.g. [G]). Let $L,L_A, L_{\bar A}$ denote a dual basis
relative to $\theta,\theta^{A},\theta^{\bar A}$. Thus, the $L_{\bar A}$
form a basis for the CR vector fields, $L_A=\overline{L_{\bar A}}$,
and $L_n$ spans the Levi nullbundle $\frak N$. We shall use the
notation introduced in \S 1 for the second (the Levi form) and third
order tensors (relative to the bases chosen). Since $M$ is Levi
uniform of rank $n-1$ near $p_0$, the Levi form $(g_{\bar A B})$
satisfies \thetag{1.9} with $r={n-1}$, 
and the third order tensor $(h_{\bar A\bar B n})$
satisfies, by Proposition 1.11, $h_{\bar n \bar B n}=h_{\bar A\bar n
n}=0$ in a neighborhood of $p_0$. Moreover, the matrix
$(h_{\bar\alpha\bar \beta n})$ is symmetric.
$(g_{\bar\alpha\beta})$ is an invertible Hermitian matrix and we shall
denote its inverse by $(g^{\alpha\beta})$,
i.e. $g_{\bar\alpha\beta}g^{\bar\alpha \gamma}=\delta^{\gamma}_\beta$.
Using the formal integrability of $\Cal V$  and
the fact that $\theta$ is real, we obtain
$$
d\theta=ig_{\bar
\alpha\beta}\theta^{\bar
\alpha}\wedge\theta^{\beta}+\phi\wedge\theta,\tag
2.9
$$
where $\phi$ is a real 1-form.
Similarly by the formal integrability, we can write
$$
d\theta^\alpha=\eta^\alpha\wedge
\theta+\eta^\alpha_\beta\wedge\theta^\beta+h^\alpha_{\bar
B}\theta^{\bar B}\wedge \theta^n,\tag 2.11
$$
for some 1-forms $\eta^\alpha$ and $\eta^\alpha_\beta$, and some
matrix $(h^\alpha_{\bar B})$.

\proclaim{Lemma 2.12} In the notation introduced above, we have
$$
h^\beta_{\bar n}=0,\quad 2i g^{\bar\alpha\beta}h_{\bar\alpha\bar\gamma
n}=h^\beta_{\bar\gamma}.\tag 2.13
$$
\endproclaim
\demo{Proof} Recall the operators $\Cal T_{\bar A}$ introduced in \S
1. Observe that, by the definition of $T''M$ and the fact
that $L_n$ spans $\frak N$, the 1-forms
$\theta$ and $\Cal T_{\bar\alpha}\theta$ form a basis for the sections of
$T''M$, and $\Cal T_{\bar n}\theta=c \theta$ for some
smooth function 
$c$. Indeed, we have
$$
\theta^\beta=g^{\bar\alpha\beta} \Cal
T_{\bar\alpha}\theta+c^\beta\theta,\tag 2.14
$$
for some smooth functions $c^\beta$, as is straightforward to
verify. 
In view of \thetag{2.14}, a direct calculation using the
definition of the third order tensor shows
$$
\left<L_{\bar B}\lrcorner d\theta^\beta, L_n\right> =
2ig^{\bar\alpha\beta} h_{\bar\alpha\bar B n}.\tag 2.15
$$
On the other hand, a calculation using \thetag{2.11} shows
$$
\left<L_{\bar B}\lrcorner d\theta^\beta, L_n\right> =
h^\beta_{\bar B},\tag 2.16
$$
which completes the proof.\qed\enddemo

Thus, we rewrite \thetag{2.11} as follows
$$
d\theta^\alpha=\eta^\alpha\wedge
\theta+\eta^\alpha_\beta\wedge\theta^\beta+h^\alpha_{\bar
\beta}\theta^{\bar \beta}\wedge \theta^n,\quad
2ig^{\bar\nu \alpha}h_{\bar\nu\bar\beta
n}=h^\alpha_{\bar\beta}. \tag 2.17
$$
In view of \thetag{2.8}, we can write
$$
d\pmatrix\omega\\ \omega^\alpha\\ \omega^{\bar \alpha}
\endpmatrix=\pmatrix \Delta&0&0\\\Delta^\alpha&\Delta^\alpha_\beta&0
\\\Delta^{\bar \alpha}&0&\Delta^{\bar \alpha}_{\bar\beta}\endpmatrix
\wedge\pmatrix \omega\\ \omega^\beta\\
 \omega^{\bar \beta}\endpmatrix+\pmatrix i\hat
g_{\bar\mu\nu}\omega^{\bar\mu}\wedge\omega^{\nu}\\ \hat h^\alpha_{\bar
\mu}\omega^{\bar \mu}\wedge \theta^n\\\overline{\hat h^\alpha_{\bar
\mu}}\omega^\mu\wedge \theta^{\bar n}\endpmatrix,\tag 2.18
$$
where $\Delta$, $\Delta^\alpha$, $\Delta^{\bar
\alpha}:=\overline{\Delta^\alpha}$, $\Delta^\alpha_\beta$,
$\Delta^{\bar \alpha}_{\bar\beta}:=\overline{\Delta^\alpha_\beta}$ are
Maurer--Cartan forms for $G$ modulo
$\omega,\omega^\alpha,\omega^{\bar\alpha}, \theta^n,\theta^{\bar n}$,
and where $\hat g_{\bar\alpha\beta}$ and $\hat h^\alpha_{\bar\beta}$
are functions on $Y\cong M\times G$ satisfying the following
$$
u^{-1}\hat
g_{\bar\alpha\beta}(p,ST)\overline{u^{\alpha}_{\mu}}u^{\beta}_\nu=\hat
g_{\bar \mu\nu}(p,T),\quad \hat
h^\alpha_{\bar\beta}(p,ST)\overline{u^{\beta}_{\mu}} =\hat
h^\nu_{\bar\mu}(p,T)u^{\alpha}_\nu,\tag 2.19
$$
where $S, T\in G$ with $S$ given by \thetag{2.6} and $p\in M$. The
last action of $G$ in \thetag{2.19} can also be described using the
third order tensor $\hat h_{\bar\alpha\bar\beta
n}=\hat g_{\bar\alpha\mu}h^\mu_{\bar\alpha}/2i$ by
$$
u^{-1}\hat h_{\bar\alpha\bar\beta n}(p,ST)
\overline{u^{\alpha}_{\mu}}\overline{u^{\beta}_{\nu}}=\hat
h_{\bar\mu\bar\nu
n}(p,T).\tag 2.20
$$
In what follows, we shall simply denote $\hat h_{\bar\alpha\bar\beta n}$
by $\hat h_{\bar \alpha\bar \beta}$.
We shall now proceed under two assumptions. The first is the following.

\definition{Condition 2.21} The matrix $(g_{\bar\alpha\beta}(p_0))$ is
{\it definite}. \enddefinition

In view of the first identity in \thetag{2.19}, we can make an
initial choice of the $\theta,\theta^{\alpha}$ so that the matrix
$(g_{\bar\alpha\beta})$ is constant with
$g_{\bar\alpha\beta}=\delta_{\bar\alpha\beta}$. Denote by $G'$ the
subgroup of $G$
consisting of those matrices $S$, as in \thetag{2.6}, for which
$$
u^{-1}
g_{\bar\alpha\beta}\overline{u^{\alpha}_{\mu}}u^{\beta}_\nu=
g_{\bar \mu\nu},\tag  2.22
$$
i.e. those for which $u>0$ and $(u^{-1/2} u^{\alpha_\beta})$ is
unitary. As a consequence of E. Cartan's work on Lie groups (see also
[S] for an elementary proof of precisely the following statement), the
(complex)
symmetric matrix $(h_{\bar\alpha\bar\beta})$ can be conjugated by
the action $v^{\alpha}_\beta\mapsto h_{\bar\alpha\bar\beta
}\overline{v^{\alpha}_\mu}\overline{v^\beta_\nu}$, where
$(v^\alpha_\beta)$ is unitary, to
the form 
$$
D_{n-1}(\lambda_{\bar 1},\ldots,\lambda_{\bar n-\bar 1}):=
\pmatrix \lambda_{\bar 1} & 0        &\ldots& 0\\
          0        & \lambda_{\bar 2}  &\ldots& 0\\
          \vdots   & \vdots         &\ddots& \vdots\\
          0        & 0
&\ldots&\lambda_{\bar n-\bar 1}\endpmatrix,\tag 2.23 $$ where
$\lambda_{\bar 1}\geq\lambda_{\bar 2}\geq\ldots\geq \lambda_{\bar
n-\bar 1}\geq 0$ are uniquely determined by
$(h_{\bar\alpha\bar\beta })$. Moreover, the subgroup of the unitary
matrices $(v^{\alpha_\beta})$ for which $$ \hat
g_{\bar\alpha\beta}\overline{v^{\alpha}_{\mu}}v^{\beta}_\nu=\hat
g_{\bar \mu\nu},\quad \hat
h_{\bar\alpha\bar\beta }\overline{v^{\alpha}_{\mu}}
\overline{v^{\beta}_{\nu}}=\hat h_{\bar\mu\bar \nu }\tag 2.24 $$
can be described explicitely (see e.g. [E4, Lemma 5.24]). Our
second assumption is the following.

\definition{Condition 2.25} The positive numbers
$\lambda_{\bar 1}(p_0),\ldots,\lambda_{\bar n-\bar 1}(p_0)$
associated with the matrix $(h_{\bar\alpha\bar \beta }(p_0))$, as
described above, are all distinct and nonzero, i.e. $$
\lambda_{\bar 1}(p_0)>\ldots>\lambda_{\bar n-\bar 1}(p_0)>0.\tag
2.26 $$ \enddefinition

Under this assumption, there is a small neighborhood of $p_0$ for
which $\lambda_{\bar 1}(p)>\ldots>\lambda_{\bar n-\bar 1}(p)>0$ for
$p$ in that
neighborhood. In view of the above and using the fact that
$\hat h_{\bar\alpha\bar\beta}(p)$ can be changed by a scalar multiple by
changing $\theta^n(p)$ by a scalar multiple, we see that is possible to
choose the basis
$\theta,\theta^{\alpha},\theta^n$ so that $g_{\bar\alpha\beta}$
is constant, equal to $\delta_{\bar\alpha\beta}$, and
$h_{\bar\alpha\bar\beta }$ satisfies
$$
\tr \hat h_{\bar\alpha\bar\beta }(p)=1,\tag
2.27
$$
for $p$ in a
neighborhood of $p_0$ (which we from now on identify with $M$). Here,
$\tr$ denotes the usual trace of a matrix; thus, we have $\tr\hat
h_{\bar\alpha\bar\beta }=\lambda_{\bar 1}+\ldots+\lambda_{\bar n-\bar
1}$.
It follows from [E4, Lemma 5.24] that the only unitary matrices
$(v^\alpha_\beta)=(u^{-1/2} 
u^{\alpha}_\beta)$
now satisfying \thetag{2.24} are the diagonal matrices
$D_{n-1}(\epsilon_1,\ldots, \epsilon_{n-1})$, where we have used the
notation introduced in \thetag{2.23} and each
$\epsilon_j\in\{-1,1\}$. Let us denote by 
$G'_1$ the subgroup of $G'$ consisting of those $S$, as in
\thetag{2.6}, where $(u^{-1/2} u^{\alpha\beta})$ equals 
$D_{n-1}(\epsilon_1,\ldots, \epsilon_{n-1})$ with
$\epsilon_j\in\{-1,1\}$. Observe that the covectors
$\omega^\alpha(p_0)$ define an ordered set of linear forms on
$\bar\Cal V_{p_0}$, each of which is invariantly defined up to
multiplication by $\pm u^{-1/2}$ 
under the action of the group
$G'_1$ in \thetag{2.5} . We
shall refer to a choice of each of these linear forms (up to
multiplication by a {\it positive} real number) as a
choice of {\it orientation for the 
normalized CR structure} at $p_0$. Thus, given such a choice of
orientation, the 
only matrices $S\in 
G'_1$ which preserve this orientation are the ones for which
$(u^{-1/2} u^{\alpha\beta})$ equals the identity. In what follows, we
shall assume that such a choice of orientation has been made. Let us
also remark that if $n=2$,
then there is only one linear form $\omega^1(p_0)$ and its annihilator
coincides with the 
Levi nullspace $\frak N_{p_0}$. Hence, in this case, normalizing the
second and third order tensor $\hat g_{\bar\alpha\beta}$, $\hat
h_{\bar\alpha\bar\beta}$ as above determines, up to multiplication by
a nonzero real
number, a linear form defining
the Levi nullspace $\frak N_{p_0}$. A choice of orientation for the
normalized CR structure at 
$p_0$ determines the sign of this real number. For this reason we
shall refer to a choice of orientation for the normalized CR structure
in the case $n=2$ 
as a choice of {\it orientation for the Levi nullspace}.

We denote by $G_1$ the subgroup of $G'_1\subset G'$ consisting of
those $S$, as in 
\thetag{2.6}, where $(u^{-1/2} u^{\alpha\beta})$ equals the
identity. It is easy to compute the Lie algebra $\frak g_1$ of
$G_1$. We have $T\in\frak g_1$ if and only if
$$
T=\pmatrix 2v&0&0\\
v^\alpha & v\delta^\alpha_\beta&0 \\ \overline{v^\alpha} &0&
v\delta^\alpha_\beta
\endpmatrix,\quad (v,v^\alpha)\in \bR
\times\bC^{n-1}.\tag 2.28
$$
Also, denote by $Y_1\subset Y$ the principal bundle over $M$ with
group $G_1$ consisting of those
$(\omega,\omega^\alpha,\omega^{\bar\alpha})^\tau$ for which $\hat
g_{\bar\alpha\beta}$ and $\hat h^{\alpha}_{\bar\mu}$ in the
structure equation \thetag{2.18} satisfy
$$
\hat g_{\bar\alpha\beta}(p,S)=\delta_{\bar\alpha\beta},\quad
\hat h_{\bar\alpha\bar
\beta }=\lambda_{\bar \beta}\delta_{\bar \alpha\bar\beta},\tag 2.29
$$
i.e.
$Y_1\cong M\times G_1\subset M\times G\cong Y$
where the trivializations are the ones obtained by choosing the basis
$\theta,\theta^\alpha, \theta^n$
as described right before \thetag{2.27}. $Y_1$ is called a reduction
of $Y_2$ (cf. e.g. [S, Chapter VII]).
In what follows, we shall,
in order to make the formulas more invariant, continue to use the
notation $g_{\bar\alpha\beta}$, $\hat h_{\bar\alpha\bar\beta }$
rather than the special form \thetag{2.29}.

Taking the pullbacks of all the forms
$\omega,\omega^\alpha,\omega^{\bar\alpha}, \theta^n,\theta^{\bar n},
\Delta, \Delta^\alpha,\Delta^{\alpha}_{\beta}$
to $Y_1$, we see from formula
\thetag{2.8} (with $S\in G_1$ now) that, pulled back to $Y_1$, the matrix
$$
\underline{\Delta}=
\pmatrix
\Delta&0&0\\\Delta^\alpha&\Delta^\alpha_\beta&0
\\\Delta^{\bar \alpha}&0&\Delta^{\bar \alpha}_{\bar \beta}\endpmatrix
$$
is valued in $\frak g_1$ modulo
$(\omega,\omega^\alpha,\omega^{\bar\alpha}, \theta^n,\theta^{\bar
n})$, i.e.
$$
\Delta^\alpha_\beta=\Delta^{\bar\alpha}_{\bar\beta}=\frac{1}{2}\Delta
\delta^\alpha_\beta\quad\mod
(\omega,\omega^\alpha,\omega^{\bar\alpha}, \theta^n,\theta^{\bar
n}).\tag 2.30
$$
Let us therefore rewrite \thetag{2.18} as follows
$$
\multline
d\pmatrix\omega\\ \omega^\alpha\\ \omega^{\bar \alpha}
\endpmatrix=\pmatrix \Delta&0&0\\\Delta^\alpha&\frac{1}{2}\Delta
\delta^\alpha_\beta&0
\\\Delta^{\bar\alpha}&0&\frac{1}{2}\Delta
\delta^\alpha_\beta\endpmatrix
\wedge\pmatrix \omega\\ \omega^\beta\\
 \omega^{\bar \beta}\endpmatrix+\pmatrix i\hat
g_{\bar\mu\nu}\omega^{\bar\mu}\wedge\omega^{\nu}\\
\hat t^\alpha_{\bar\mu\nu} \omega^{\bar\mu}\wedge\omega^\nu\\\overline{
\hat t^\alpha_{\bar\mu\nu}} \omega^{\mu}\wedge\omega^{\bar \nu}
\endpmatrix\\+ \pmatrix 0\\ \hat h^\alpha_{\bar
\mu}\omega^{\bar \mu}\wedge \theta^n+\hat
q^\alpha_{\nu\beta}\omega^\beta\wedge\omega^\nu+\hat
r^{\alpha}_\nu\omega^\nu\wedge\theta^n+\hat s^{\alpha}_
\nu\theta^{ \bar n}\wedge \omega^{\nu}\\ \overline{\hat
h^\alpha_{\bar \mu}}\omega^\mu\wedge \theta^{\bar
n}+\overline{\hat
q^\alpha_{\nu\beta}}\omega^{\bar\beta}\wedge\omega^{\bar\nu}+\overline
{\hat r^{\alpha}_\nu}\omega^{\bar\nu}\wedge\theta^{\bar n}+
\overline{\hat s^{\alpha}_\nu} \theta^{n}\wedge \omega^{\bar
\nu}\endpmatrix,\endmultline\tag 2.31 $$ where, by a slight abuse
of notation, all the forms in \thetag{2.31} denote the pulled
back forms on $Y_1$; we require the $\hat q^\alpha_{\nu\beta}$ to
be skew symmetric in $\nu$, $\beta$. The 1-form $\Delta$ is
uniquely determined by the form of the structure equation
\thetag{2.31} up to transformations $$ \tilde
\Delta=\Delta+a\omega,\tag 2.32 $$ where $a$ is a smooth function
on $Y_1$; i.e. replacing $\Delta$ by $\tilde \Delta$ given by
\thetag{2.32} preserves \thetag{2.31} and this is the only
transformation preserving \thetag{2.31}, as is easily verified
(cf. also [G, Lecture 3]). In fact, a direct calculation shows
that $$ \Delta=u^{-1}du+i\hat
g_{\bar\mu\nu}u^{-1/2}(\overline{u^\mu}\omega^\nu-u^\nu\omega^{\bar
\mu})+\phi\mod\omega. $$ where $\phi$ is the pullback of a real
1-form on $M$. The 1-form $\theta^n$ on $Y_1$ is determined by the
condition $\tr \hat h_{\bar\alpha\bar\beta}\equiv 1$ up to
transformations $$ \tilde
\theta^n=\theta^n+c_\beta\omega^\beta+c\omega,\tag 2.33 $$ where
$c,c_\beta$ are smooth functions on $Y_1$. The 1-forms
$\Delta^\alpha$ are determined modulo
$\omega,\omega^\nu,\omega^{\bar \mu},\theta^n,\theta^{\bar n}$.
The precise form of the indeterminacy in $\Delta^\nu$ is not important
at this point.

By using the integrability of $\Cal V$ and the fact that $\tr \hat
h_{\bar\alpha\bar \beta}$ is constant on the fibers $Y_1\to M$, we
deduce that we can write $$ \aligned d\tilde\theta^n &=
d\theta^n+dc_\beta\wedge\omega^\beta+c_\beta
d\omega^\beta+dc\wedge\omega +cd\omega\\
&=\eta\wedge\omega+\eta_\beta\wedge\omega^\beta+\tilde k_{\bar
\mu}\, \omega^{\bar \mu}\wedge\tilde\theta^n+\hat k\,
\tilde\theta^{\bar n}\wedge\tilde \theta^n,\endaligned\tag 2.34
$$ where $\eta$ and $\eta_\beta$ are 1-forms (depending also on
$c_\beta$ and $c$),  $$ \tilde k_{\bar\mu} =\hat k_{\bar\mu}
+c_\nu\hat h^\nu_{\bar\mu}-\overline{c_{\mu}}\hat k,\tag 2.35$$
and $\hat k_{\bar\mu}$, $\hat k$ are smooth functions on $Y_1$.
Moreover, $\hat k$ is constant on the fibers $Y_1\to M$, whereas
$\hat k_{\mu}(p,ST)=u^{-1/2} \hat k_{\bar\mu}(p,T)$ for $S,T\in
G_1$ and $S$ of the form \thetag{2.6}. The reader can easily
verify that the function $\hat k$ is also uniquely determined,
i.e. independent of the choice of $\Delta,\Delta^\nu,\theta^n$. It
is also easy to check that the functions $\hat r^\alpha_\nu$,
$\hat s^\alpha_\nu$ in \thetag{2.31} are uniquely determined, and
that they are constant on the fibers $Y_1\to M$. Hence, $\hat
k(p)$, $\hat r^\alpha_\nu(p)$, $\hat s^\alpha_\nu(p)$ are
invariants of the CR structure on $M$ at $p\in M$. However, they
are not independent as the following proposition shows.

\proclaim{Proposition 2.36} We have the following identities $$
\hat g_{\bar\mu \eta} \hat s^\eta_\nu= \hat g_{\bar\xi \nu}
\overline{ \hat r^\xi_\mu}=\overline{\hat g_{\bar\nu \xi} \hat
r^\xi_\mu},\tag 2.37 $$ and $$ \hat k=\tr\left ( \hat
h_{\bar\alpha\mu}\overline{\hat r^\mu_\beta}+\hat
h_{\bar\mu\bar\beta} \overline{\hat
r^\mu_\alpha}\right)_{\alpha\beta}.\tag 2.38 $$
\endproclaim
\demo{Proof} Observe that the second identity in \thetag{2.37}
follows directly from the fact that $(\hat g_{\bar\alpha\beta})$
is Hermitian. Differentiating the first row in the structure
equation \thetag{2.31}, using the facts that $d^2\omega=0$ and
$\hat g_{\bar\alpha\beta}$ is constant on $Y_1$, we obtain $$
0=d\Delta\wedge\omega-\Delta\wedge d\omega+i\hat
g_{\bar\mu\nu}(d\omega^{\bar\mu}\wedge
\omega^\nu-\omega^{\bar\mu}\wedge d\omega^{\nu}).\tag 2.39 $$
Applying equation \thetag{2.31} again, we obtain $$ \multline
0=\left(d\Delta+i\hat g_{\bar\mu\nu}(\Delta^\nu\wedge
\omega^{\bar\mu}-\Delta^{\bar\mu}\wedge\omega^\nu)\right)\wedge\omega\\
+i\hat g_{\bar\mu\nu}(\overline{\hat
t^\mu_{\bar\xi\eta}}\omega^\xi\wedge\omega ^{\bar\eta}
\wedge\omega^\nu-\hat t^\nu_{\bar\xi\eta}\omega^ {\bar
\xi}\wedge\omega^{\eta} \wedge\omega^{\bar\mu})  +i\left(\hat
g_{\bar\xi\nu}\overline{\hat s^\xi_\mu}-\hat g_{\bar\mu\eta} \hat
r^\eta_\nu\right)\theta^n\wedge
\omega^{\mu}\wedge\omega^\nu\\+i\left(\hat g_{\bar\mu\eta} \hat
s^\eta_\nu-\hat g_{\bar\xi\nu}\overline{\hat
r^{\xi}_\mu}\right)\theta^{\bar n} \wedge
\omega^{\bar\mu}\wedge\omega^\nu  - i\hat
g_{\bar\mu\nu}(\overline{\hat h^{\mu}_{\bar \xi}}\theta^{\bar
n}\wedge\omega^{\xi}\wedge\omega^\nu + \hat h^{\nu}_{\bar
\eta}\theta^{n}\wedge\omega^{\bar \mu}\wedge\omega^{\bar
\eta}).\endmultline \tag 2.40 $$ The first identity in
\thetag{2.37} follows immediately from \thetag{2.40}. To prove
\thetag{2.38}, we differentiate the formula for $d\omega^\alpha$
given by \thetag{2.31}. We obtain $$\multline
0=d\Delta^\alpha\wedge\omega-\Delta^\alpha\wedge
d\omega+\frac{1}{2}(d \Delta\wedge\omega^\alpha-\Delta\wedge
d\omega^\alpha)+ d\hat
t^\alpha_{\bar\mu\nu}\wedge\omega^{\bar\mu}\wedge\omega^\nu+\\
\hat
t^\alpha_{\bar\mu\nu}(d\omega^{\bar\mu}\wedge\omega^\nu-\omega^{\bar\mu}
\wedge d \omega^\nu)+ d\hat
h^\alpha_{\bar\mu}\wedge\omega^{\bar\mu}\wedge\theta^n+ \hat
h^\alpha_{\bar\mu}(d\omega^{\bar\mu}\wedge\theta^n-
\omega^{\bar\mu}\wedge d\theta^n)+\\ d\hat
r^\alpha_{\nu}\wedge\omega^{\nu}\wedge\theta^n+ \hat
r^\alpha_{\nu}(d\omega^{\nu}\wedge\theta^n-\omega^{\nu}\wedge
d\theta^n) +d\hat s^\alpha_{\nu}\wedge\theta^{\bar
n}\wedge\omega^{\nu}+ \hat s^\alpha_{\nu}(d\theta^{\bar
n}\wedge\omega^{\nu}-\theta^{\bar n}\wedge
d\omega^{\nu}).\endmultline\tag 2.41 $$ Let us write $d\hat
h^\alpha_{\bar\mu}=e^\alpha_{\bar\mu}\theta^{\bar n}$ modulo
$\omega,\omega^\nu,\omega^{\bar \nu},\theta^n$. If we substitute
\thetag{2.31} and \thetag{2.35} in \thetag{2.41} and collect
the $\omega^{\bar\mu}\wedge\theta^{\bar n}\wedge\theta^n$-terms,
we obtain $$ \hat k\hat h^\alpha_{\bar\mu}=\hat
h^\alpha_{\bar\gamma}\overline{\hat r^\gamma_\mu} + \hat
s^\alpha_\gamma\hat h^\gamma_{\bar \mu}-e^\alpha_{\bar\mu}.\tag
2.42 $$ Multiplying \thetag{2.42} by $2i\hat g_{\bar\beta\alpha}$,
summing over $\alpha$, and using \thetag{2.37}, we obtain $$ \hat
k\hat h_{\bar \beta\bar\mu}=\hat
h_{\bar\beta\bar\gamma}\overline{\hat r^\gamma _\mu}+\hat
h_{\bar\gamma\bar\mu}\overline{\hat r^\gamma _\beta}-2i\hat
g_{\bar\beta\alpha} e^\alpha_{\bar\mu}.\tag 2.43 $$ The identity
\thetag{2.38} now follows by taking the trace of \thetag{2.43}
and using the fact that $\tr \hat h_{\bar\alpha\bar\beta}\equiv
1$. (In particular, $\tr\hat g_{\bar\beta\alpha}
e^\alpha_{\bar\mu}\equiv 0$.)   \qed\enddemo

For future reference, let us remark that the functions $\hat
t^\alpha_{\bar\mu\nu}$ in \thetag{2.31} depend on the choice of
$\theta^n$. For a fixed such choice, the functions $\hat
t^\alpha_{\bar\mu\nu}$ satisfy the following on $Y_1$ $$ \hat
t^\alpha_{\bar\mu\nu}(p,ST)= \frac{\hat
t^\alpha_{\bar\mu\nu}(p,T)}{\sqrt{u}}+i\hat
g_{\bar\mu\nu}(p)\frac{u^\alpha}{u} +i\frac{1}{2}\hat
g_{\bar\mu\gamma}(p)\frac{u^\gamma}{u}\delta^\alpha_{\nu},\tag
2.44 $$ where $S,T\in G_1$ and $S$ is given by \thetag{2.6}. Let
us also observe that replacing $\theta^n$ by $\tilde\theta^n$, as
given by \thetag{2.33}, the $\hat t^\alpha_{\bar\mu\nu}$ change
by $$ \hat t^\alpha_{\bar\mu\nu}\to \hat t^\alpha_{\bar
\mu\nu}-(\hat h^\alpha_{\bar\mu} c_\nu+\hat
s^\alpha_\nu\overline{c_\mu}).\tag 2.45 $$

\heading 3. The 5 dimensional case\endheading

We now restrict our attention to the case $n=2$. The conditions
2.21 and 2.25 reduce to requiring that $M$ is Levi uniform of
rank one and $2$-nondegenerate at $p_0$. Moreover, we have fixed an
orientation of the Levi nullspace $\frak N_{p_0}$, as explained in section
2, in order to distinguish a component $Y_1$ of the disconnected
bundle $Y'_1$. Our aim is to define a 
uniquely determined 
submanifold $Y_2\subset Y_1$ which can be viewed as a principal fiber
bundle over $M$ whose group is a subgroup of $G_1$, and on the
principal bundle $Y_2\to M$ 
determine the forms $\theta^n$, $\Delta$, and $\Delta^1$ uniquely.
(Recall that $n=2$ so that Greek indices $\alpha$, $\beta$ run
over the single integer $1$.) We have to distinguish two cases. We
shall begin by handling the most difficult case, which does not
appear to be the generic one but which contains the tube over the
light cone as given by Example 1.7. Observe in what follows that
the condition $\tr \hat h_{\bar\alpha\bar\beta}\equiv 1$ reduces
to $\hat h_{\bar 1\bar 1}\equiv 1$. To simplify the notation, we
shall use the fact that $\hat g_{\bar 1 1}\equiv 1$. Also, recall the
invariant $\hat k$ defined by \thetag{2.34}.

\subhead 3.1. Case 1: $|\hat k(p_0)|=2$\endsubhead Consider the
formula in \thetag{2.35} for the coefficient $\tilde
k_{\bar 1}$ of $\omega^{\bar 1}\wedge\tilde \theta^2$ in $d\tilde
\theta^2$ and rewrite it in the following form $$\tilde k_{\bar
1}=\hat k_{\bar 1}+2ic_1-\overline{c_1}\hat k.\tag 3.1.1 $$ Since
$|\hat k(p_0)|=2$, it is not possible to solve the equation
$\tilde k_{\bar 1}=0$ for $c_1$ in a neighborhood of $p_0$. Let us
write $$\hat k=2ire^{i 
t},\tag 3.1.2 $$ where $r$ and $t$ are real valued functions on
$Y_1$ which are constant on the fibers of $Y_1\to M$. We shall seek
$c_1$ in the form $$ c_1=e^{it/2}(\rho_1+i\rho_2),\tag 3.1.3 $$
where $\rho_1,\rho_2\in\bR$. Since $r(p_0)=1$ we can choose
$\rho_2$ uniquely so that $$ \tilde k_{\bar 1}= ir'e^{it/2}, $$ for
some real valued function $r'$ on $Y_1$ which is
constant on the fibers of $Y_1\to M$. The function $\rho_2$
satisfies the same transformation rule as $\hat k_{\bar 1}$, i.e.
$$\rho_2(p,ST)=u^{-1/2}\rho_2(p,T),\tag 3.1.4$$ for $S,T\in G_1$
and $S$ of the form \thetag{2.6}. We can express the above by
saying that $\im e^{-it/2} c_1$ is uniquely determined by the
condition 
$$ 
\re e^{-it/2}\hat k_{\bar 1}\equiv 0.\tag 3.1.5 
$$
Now, with $\im e^{-it/2} c_1$ determined by \thetag{3.1.5}, the
equation \thetag{2.45} implies that the function $\hat t^1_{\bar 11}$,
given by \thetag{2.31}, is determined up to $$\hat t^1_{\bar
11}\to \hat t^1_{\bar 11} -\left(2i\rho_1 e^{it/2}+\hat
s^1_1 \rho_1 e^{-it/2}\right) $$ or equivalently, in view of
Proposition 2.36 and \thetag{3.1.2}, $$ \aligned\hat t^1_{\bar
11}\to  & \hat t^1_{\bar 11} -2i\left(\rho_1 e^{it/2}+\frac{1}{2}\,
\hat k \rho_1 e^{-it/2}\right)\\=&\hat t^1_{\bar 11} -2i\rho_1
e^{it/2}\left(1+\frac{1}{2}\, r\right).\endaligned\tag 3.1.6 $$
Let us write $u^1=e^{it/2}(x+iy)$. It follows from \thetag{2.44}
that there is a submanifold $Y_2\subset Y_1$, which is defined
(uniquely in view of \thetag{3.1.6}) by the equation $$ \re
(e^{-it/2}\hat t^1_{\bar 11})=0.\tag 3.1.7 $$ The manifold $Y_2$
can be viewed as a principal fiber bundle $Y_2\to M$ with group
$G_2\subset G_1$, where $G_2$ is defined by the equation $$ \im
e^{-it(p_0)/2} u^1=0,\tag 3.1.8 $$ if we let $G_2$ act on $Y_2$
as follows $$ g\pmatrix
\omega\\\omega^\alpha\\\omega^{\bar\alpha}\endpmatrix:= \pmatrix
u&0&0\\ e^{it/2}x&\sqrt{u}&0\\e^{-it/2}x&0&\sqrt{u}\endpmatrix
\pmatrix \omega\\\omega^\alpha\\\omega^{\bar\alpha}\endpmatrix,
\tag 3.1.9 $$ for $g\in G_2\subset G_1$ of the form $$
\pmatrix u&0&0\\
e^{it(p_0)/2}x&\sqrt{u}&0\\e^{-it(p_0)/2}x&0&\sqrt{u}\endpmatrix.\tag
3.1.10 $$ Observe that the principal $G_2$ bundle $Y_2$ is a
reduction of the principal $G_1$ bundle $Y_1$ only if the phase
function $t$ is constant. By definition, we have $\hat t^1_{\bar
11}=ir''e^{it/2}$, for some real valued function $r''$, on the
bundle $Y_2$. Hence, we can determine $\rho_1$ uniquely in
\thetag{3.1.6} so that $$ \hat t^1_{\bar 11}\equiv 0,\tag 3.1.11
$$ on $Y_2$.  It follows from \thetag{2.44} and \thetag{3.1.6}
that we have
$$\rho_1(p,ST)=\frac{\rho_1(p,T)}{\sqrt{u}}+\frac{3}{2(2+r)}\,\frac{x}{u},
\tag 3.1.12$$ for $S,T\in G_2$ and $S$ of the form \thetag{2.6}
with $u^1=e^{it/2}x$ and $x\in \bR$.

Thus, $\theta^2$ is determined on $Y_2$ up to 
$$
\tilde\theta^2=\theta^2+c\omega\tag 3.1.13 
$$ 
by \thetag{3.1.5} and 
\thetag{3.1.11}. Observe that on $Y_2$, where $u^1=e^{it/2}x$ for
$x\in\bR$, we have $$ \Delta^1=e^{it/2}\xi\quad\mod
\omega,\omega^1,\omega^{\bar 1},\theta^2, \theta^{\bar 2},\tag
3.1.14 $$ where $\xi$ is a real 1-form on $Y_2$ which is not
uniquely determined. We can write the structure equation for
$d\omega^1$ on $Y_2$ as follows $$\multline d\omega^1=
e^{it/2}\xi\wedge\omega+\frac{1}{2}\Delta\wedge\omega^1+\hat
h^1_1\omega^{\bar 1}\wedge\theta^2+\hat
r^1_1\omega^1\wedge\theta^2+\hat s^1_1\theta^{\bar
2}\wedge\omega^1\\+ie^{it/2}(b\omega^1\wedge\omega+\bar
b\omega^{\bar 1}\wedge\omega+e\theta^2\wedge\omega+\bar
e\theta^{\bar 2}\wedge\omega),\endmultline\tag 3.1.15 $$ for some
functions $b$ and $e$ on $Y_2$. Using \thetag{2.32} and
\thetag{3.1.13}, we obtain $$\multline d\omega^1= e^{it/2}
\tilde \xi\wedge\omega+\frac{1}{2}\tilde \Delta\wedge\omega^1+\hat
h^1_1\omega^{\bar 1}\wedge\tilde\theta^2+\hat
r^1_1\omega^1\wedge\tilde\theta^2+\hat s^1_1\tilde\theta^{\bar
2}\wedge\omega^1\\+ie^{it/2}(\tilde
b\omega^1\wedge\omega+\overline{\tilde b }\omega^{\bar
1}\wedge\omega+ e\tilde\theta^2\wedge\omega+\bar
e\tilde\theta^{\bar 2}\wedge\omega),\endmultline\tag 3.1.16 $$
where 
$$
\tilde b=b+\frac{1}{2}\left(\frac{1}{2} ae^{-it/2}+\bar
c\hat s^1_1 e^{-it/2}+ \bar c\overline{\hat h^1_1} e^{it/2}- c\hat
r^1_1 e^{-it/2}\right)\tag 3.1.17
$$ 
and 
$$
\multline\tilde
\xi=\xi+\frac{1}{2}\left(\frac{1}{2}( ae^{-it/2}\omega^1+\bar
ae^{it/2}\omega^{\bar 1})+\bar c\hat s^1_1 e^{-it/2}\omega^ 1 +
c\overline{\hat s^1_1} e^{it/2}\omega^{\bar 1}\right)\\-
\frac{1}{2}\left( \bar c\overline{\hat h^1_1} e^{it/2}\omega^1+
c\hat h^1_1 e^{-it/2}\omega^{\bar 1}+ c\hat r^1_1
e^{-it/2}\omega^1+\bar c\overline{\hat r^1_1} e^{it/2}\omega^{\bar
1}\right)+q\omega.\endmultline\tag 3.1.18
$$ 
Here, $q$ is an
arbitrary real valued function on $Y_2$, and $a$ and $c$ are as in
\thetag{2.32} and \thetag{3.1.13} respectively. We deduce that
$\xi$ is determined by the structure equation \thetag{3.1.16} up
to transformations given by \thetag{3.1.18}. Let us rewrite
\thetag{3.1.17} and \thetag{3.1.18} using Proposition 2.36 and
\thetag{3.1.2} as follows $$\tilde
b=b+\frac{1}{2}e^{-it/2}\left(\frac{1}{2}a+i\left(r-2\right)
\bar \zeta+ir\zeta\right),\tag 3.1.19 $$ and
$$\multline\tilde \xi=\xi+\frac{1}{2}e^{-it/2}\left(\frac{1}{2}
a+i\left(r+2\right)\bar \zeta + ir\zeta
\right)\omega^1\\+\frac{1}{2}e^{it/2}\left(\frac{1}{2}
a-i\left(r+2\right) \zeta - ir\bar \zeta
\right)\omega^{\bar 1}+q\omega,\endmultline\tag 3.1.20$$ where we
have used the notation  $$\zeta=e^{-it}c.$$ By \thetag{2.40}, we
have on $Y_1$, $$\aligned d\Delta &=i\hat g_{\bar 11}(\Delta^{\bar
1}\wedge\omega^1-\Delta^1\wedge\omega^{\bar
1})+\Phi\wedge\omega\\&=i(\Delta^{\bar
1}\wedge\omega^1-\Delta^1\wedge\omega^{\bar
1})+\Phi\wedge\omega,\endaligned\tag 3.1.21 $$ for some real
1-form $\Phi$; on the second line of \thetag{3.1.21}, we have
used $\hat g_{\bar 11}=1$. Hence, on $Y_2$ we have, by
\thetag{3.1.14}, $$d\Delta=ie^{-it/2}\xi\wedge\omega^1-i
e^{it/2}\xi\wedge\omega^{\bar 1} +\Psi\wedge\omega
+\Psi_1\wedge\omega^1+\Psi_{\bar 1}\wedge\omega^{\bar 1},\tag
3.1.22$$ for some 1-forms $\Psi,\Psi_1$, $\Psi_{\bar
1}=\overline{\Psi_1}$ on $Y_2$ such that $\Psi_1=0$ modulo
$\omega^1,\omega^{\bar 1},\theta^2,\theta^{\bar 2}$. Recall that
$\Delta$ is determined up to transformations \thetag{2.32},
$\theta^2$ up to transformations \thetag{3.1.13}, and $\xi$ up to
transformations \thetag{3.1.20}. Substituting in
\thetag{3.1.22}, we obtain $$\aligned d\tilde\Delta&
=d\Delta+da\wedge\omega+ad\omega\\ &=i e^{-it/2}\tilde
\xi\wedge\omega^1 -i e^{it/2}\tilde\xi\wedge\omega^{\bar 1} +i\tilde
f\omega^{\bar 1}\wedge\omega^1+\ldots,\endaligned\tag 3.1.23$$
where $\ldots$ signify the remaining terms in the expansion of
$d\tilde\Delta$ and $$\tilde f=f-\frac{1}{2}a+i\zeta-i\bar
\zeta,\tag 3.1.24$$ for some real valued function $f$ on $Y_2$.
Hence, $a$ is uniquely determined as a function of $\zeta$ by the
condition $$\tilde f=0,\tag 3.1.25$$ and we have $$a=2i(\zeta-\bar
\zeta).\tag 3.1.26$$ In view of \thetag{3.1.3}, \thetag{3.1.4},
and \thetag{3.1.12}, we have $$ dc_1=e^{it/2}\frac{3}{2(2+r)}\xi\
\mod\Delta,\omega,\omega^1,\omega^{\bar 1}, \theta^2,\theta^{\bar
2}.\tag 3.1.27 $$ It follows, by using \thetag{2.34},
\thetag{3.1.20} and \thetag{3.1.24}, that
$$
d\tilde\theta^2=e^{it/2}\frac{3}{2(2+r)}\tilde
\xi\wedge\omega^1+\tilde m\omega^{\bar
1}\wedge\omega^1+\ldots,\tag 3.1.28
$$ 
where $\ldots$ signify the
remaining terms in the expansion of $d\tilde \theta^2$ and $\tilde
m$ is given by 
$$
\aligned \tilde m
&=m+\frac{3ie^{it}}{4}\frac{(r+1)}{(r+2)}\left(\zeta+\bar\zeta
\right)+ic\\
&=m+\frac{ie^{it}}{4(r+2)}\left((7r+11)\zeta+3(r+1)\bar\zeta\right)
\endaligned\tag 3.1.29 
$$ 
for some function $m$ on $Y_2$; in the last line of \thetag{3.1.29},
we have used $c=e^{it}\zeta$. Since $r(p_0)=1$, we can determine
$\zeta$ uniquely by the condition $$\tilde m=0.\tag 3.1.30$$

Let us summarize our efforts so far. We have determined
$\Delta$ and $\theta^2$ uniquely on $Y_2$. The 1-form $\xi$ is determined
up to $$\tilde\xi=\xi+q\omega,\tag 3.1.31$$ for some real valued
function $q$ on $Y_2$. We shall conclude the construction in this
section by defining a unique choice of $\xi$.
In view of
\thetag{3.1.22}, we have
$$d\Delta=ie^{-it/2}\tilde\xi\wedge\omega^1
-ie^{it/2}\tilde\xi\wedge\omega^{\bar 1}+\tilde l
\omega^{1}\wedge\omega+\overline{\tilde l} \omega^{\bar
1}\wedge\omega +\ldots,\tag 3.1.32$$ where $$ \tilde l
=l+ie^{-it/2}q\tag 3.1.33$$ for some function $l$
on $Y_2$. Hence, we may determine $q$ uniquely by the condition $$
\im e^{it/2}\tilde l=0.\tag 3.1.34$$ 
The uniquely determined, linearly independent 1-forms $$\omega,\re
\omega^1,\im \omega^{ 1},\re \theta^2,\im \theta^{
2},\Delta,\xi\tag 3.1.35$$ on $Y_2$ form a global coframe for
$T^* Y_2$ and, hence, define an absolute parallelism on $Y_2$.
The following result is a consequence of the construction above
(see [G] or [CM]). We use the notation
introduced previously, and also
$$
\underline{\omega}:=\pmatrix \omega\\ \re \omega^1\\\im \omega^1\\
\re \theta^2\\ \im \theta^2\\ \Delta\\ \xi\endpmatrix.\tag 3.1.36
$$

\proclaim{Theorem 3.1.37} Let $M$ be a $5$-dimensional CR manifold of
hypersurface type which is $2$-nondegenerate and Levi uniform of rank
$1$ at $p_0\in M$. Suppose that $|\hat k(p_0)|=2$ and that an
orientation is chosen for the Levi nullspace $\frak N_{p_0}$ (as
explained in \S $2$). Then,
there exists a principal fiber bundle 
$Y_2\to M$ with a two dimensional structure group $G_2\subset
GL(\bC^3)$ and a $1$-form
$\underline{\omega}$ on $Y_2$ which 
defines an isomorphism between $T_y Y_2$ and 
$\bR^{7}$ for every
$y\in Y_2$ and such that the following holds. Let $M'$ be a
$5$-dimensional CR manifold of 
hypersurface type which is $2$-nondegenerate and Levi uniform of rank
$1$ at $p'_0\in M$. Suppose that the invariant $|\hat k'(p_0)|=2$
(where corresponding 
objects for $M'$ are denoted with $'$) and 
that an orientation is chosen for the Levi nullspace $\frak
N'_{p'_0}$. Then, if there 
exists a local CR diffeomorphism $f\:(M,p_0)\to (M',p'_0)$ preserving the
oriented Levi nullspaces $\frak N_{p_0}$ and 
$\frak N'_{p'_0}$,
there exists a diffeomorphism $F\:Y_2\to Y_2'$ with $\pi'\circ
F(\pi^{-1}(p_0))=\{p_0'\}$ such that
$F^*\underline{\omega'}=\underline{\omega}$ and the following
diagram commutes
$$
\CD 
Y_2 @> F>> Y_2'\\ @V \pi VV @VV \pi' V\\ M @>> f> M'.
\endCD\tag 3.1.38
$$ 
Conversely, if there exists a diffeomorphism $F\:Y_2\to Y_2'$ with $\pi'\circ
F(\pi^{-1}(p_0))=p_0'$ such that
$F^*\underline{\omega'}=\underline{\omega}$, then there exists a CR
diffeomorphism $f\:(M,p_0)\to (M',p_0')$ preserving the
oriented Levi nullspaces $\frak N_{p_0}$ and 
$\frak N'_{p_0'}$, such that
\thetag{3.1.38} commutes. 

The $1$-form $\underline{\omega}$ is 
given by \thetag{3.1.36} and is uniquely determined
by \thetag{3.1.5}, \thetag{3.1.11}, \thetag{3.1.25},
\thetag{3.1.30}, and \thetag{3.1.34}. The group $G_2$ is defined by
\thetag{3.1.10}.
\endproclaim

\subhead 3.2. Case 2: $|\hat k(p_0)|\neq 2$\endsubhead First,
since $|\hat k(p_0)|\neq 2$, there is a small neighborhood of
$p_0$ in $M$ in which $|\hat k|\neq 2$. In this neighborhood,
which we identify with $M$ in this subsection, we can solve
uniquely for $c_1$ in the equation $$\tilde k_{\bar 1}=0,\tag
3.2.1$$ where $\tilde k_{\bar 1} $ is given by equation
\thetag{3.1.1}. Thus, $c_1$ is uniquely determined by the
condition \thetag{3.2.1}, and hence $\theta^2$ is determined up
to transformations of the form \thetag{3.1.13} . Moreover, it
follows from \thetag{2.44} that there exists a uniquely
determined submanifold $Y_3\subset Y_1$ defined by $$\hat
t^1_{\bar 11}=0.\tag 3.2.2$$ The submanifold $Y_3$ is a subbundle
(a reduction) of the principal $G_1$-bundle $Y_1$ with group
$G_3$, where $G_3$ is the subgroup of $G_1$ defined by
$$u^1=0.\tag 3.2.3$$ Thus, on $Y_3$ we have $\Delta^1=0$ modulo
$\omega,\omega^1,\omega^{\bar 1},\theta ^2,\theta^{\bar 2}$. It
follows that the structure equation for $d\omega^1$ on $Y_3$ can
be written $$\multline d\omega^1= \frac{1}{2}\tilde
\Delta\wedge\omega^1+\hat h^1_1\omega^{\bar 1}\wedge\tilde
\theta^2+\hat r^1_1\omega^1\wedge\tilde \theta^2+\hat s^1_1\tilde
\theta^{\bar 2}\wedge\omega^1\\+{}^1\tilde b
\omega^1\wedge\omega+{}^2b\tilde\omega^{\bar
1}\wedge\omega+{}^1e\theta^2\wedge\omega+{}^2 e\theta^{\bar
2}\wedge\omega,\endmultline\tag 3.2.4 $$ where $$ ^1\tilde b
:={}^1b-\hat r^1_1 c+\hat s^1_1\bar c+\frac{1}{2}a,\quad ^2\tilde
b :={}^2b-\hat h^1_1 c\tag 3.2.5$$ for some functions $^1b$,
$^2b$, $^1e$, and $^2e$ on $Y_3$. We can determine $c$, and hence
$\theta^2$, uniquely on $Y_3$ by the condition $$^2\tilde b=0.\tag
3.2.6$$ We then determine $a$, and hence $\Delta$, uniquely on
$Y_3$ by the condition $$\re ^1\tilde b=0.\tag 3.2.7$$ The uniquely
determined 1-forms $$\omega,\re \omega^1,\im \omega^{ 1},\re
\theta^2,\im \theta^{ 2},$$ on $Y_3$ form a
global coframe for $T^* Y_3$ and, hence, define an absolute
parallelism on $Y_3$.
As in 3.1, we have the following result. We use the notation
introduced above, and also
$$
\underline{\omega}:=\pmatrix \omega\\ \re \omega^1\\\im \omega^1\\
\re \theta^2\\ \im \theta^2\\ \Delta\endpmatrix\tag 3.2.8
$$

\proclaim{Theorem 3.2.9} Let $M$ be a $5$-dimensional CR manifold of
hypersurface type which is $2$-nondegenerate and Levi uniform of rank
$1$ at $p_0\in M$. Suppose that $|\hat k(p_0)|\neq 2$ and that an
orientation is chosen for the Levi 
nullspace $\frak N_{p_0}$ (as explained in \S$2$). Then, there exists
a principal fiber bundle 
$Y_3\to M$ with a one dimensional structure group $G_3\subset
GL(\bC^3)$ and a $1$-form
$\underline{\omega}$ on $Y_3$ which 
defines an isomorphism between $T_y Y_3$ and 
$\bR^{6}$ for every
$y\in Y_3$ and such that the following holds. Let $M'$ be a
$5$-dimensional CR manifold of 
hypersurface type which is $2$-nondegenerate and Levi uniform of rank
$1$ at $p'_0\in M$. Suppose that the invariant $|\hat k'(p_0)|\neq 2$
(where corresponding 
objects for $M'$ are denoted with $'$) and 
that an orientation is chosen for the Levi nullspace $\frak
N'_{p'_0}$. Then, if there
exists a local CR diffeomorphism $f\:(M,p_0)\to (M',p'_0)$ preserving the
oriented Levi nullspaces $\frak N_{p_0}$ and 
$\frak N'_{p_0'}$,
there exists a diffeomorphism $F\:Y_3\to Y_3'$ with $\pi'\circ
F(\pi^{-1}(p_0))=\{p_0'\}$ such that
$F^*\underline{\omega'}=\underline{\omega}$ and the following
diagram commutes
$$
\CD 
Y_3 @> F>> Y_3'\\ @V \pi VV @VV \pi' V\\ M @>> f> M'.
\endCD\tag 3.2.10
$$ 
Conversely, if there exists a diffeomorphism $F\:Y_3\to Y_3'$ with $\pi'\circ
F(\pi^{-1}(p_0))=p_0'$ such that
$F^*\underline{\omega'}=\underline{\omega}$, then there exists a CR
diffeomorphism $f\:(M,p_0)\to (M',p_0')$ preserving the
oriented Levi nullspaces $\frak N_{p_0}$ and
$\frak N'_{p_0'}$, such that
\thetag{3.2.10} commutes. 

The $1$-form $\underline{\omega}$ is 
given by \thetag{3.2.8} and is uniquely determined
by \thetag{3.2.1}, \thetag{3.2.6}, and \thetag{3.2.7}. The group
$G_3$ is defined by \thetag{3.2.3}, and on $Y_3$ \thetag{3.2.2}
holds. 
\endproclaim

\heading 4. A curvature characterization of the tube over the light
cone in $\bC^3$\endheading 

In this section, we shall continue to consider only the case
$n=2$. We shall change slightly the convention from previous
sections that capital Roman indices $A,B$, etc.\ run over the set
$\{1,2\}$ (i.e. $\{1,\ldots,n\}$ with $n=2$) and instead let them
run over the set $\{1,2,3\}$. We shall also use the convention
that small Roman indices $a,b,$ etc.\ run over the set
$\{0,1,2,3\}$.

Denote by $\Gamma$ the light cone in $\bR^3$, i.e. the zero locus of
the quadratic form $\{x,x\}$ where $\{\cdot,\cdot\}$ denotes the
bilinear form which in the standard coordinates of $\bR^3$ is given by
$$
\{x,y\}:=x^1y^1+x^2y^2-x^3y^3. \tag 4.1
$$
We shall denote by $\Gamma_\bC$ the tube in $\bC^3$ over $\Gamma$ (as
in Example 1.7). Hence, $\Gamma_\bC$ is the zero locus of
$\{Z,Z\}_\bC$, where $\{z,w\}_\bC$ for complex 
vectors $z, w\in\bC^3$ is defined by
$$
\{z,w\}_\bC:=\{\re z,\re w\}.\tag 4.2
$$
We shall call a {\it frame} for $\Gamma_\bC$ a 4-tuple
$(Z_0,Z_1,Z_2,Z_3)$ of 4-vectors, where
$$
Z_0=t(1,z_0),\quad Z_A=(0,x_A),\quad t\in\bR,\ z_0\in\bC^3,\
x_A\in\bR^3,\tag 4.3
$$
which satisfy the following conditions. The real
vectors $x_1,x_2,x_3\in\bR^3$ satisfy
$$
\aligned
\{x_1,x_1\} &=\{x_3,x_3\}=\{x_1,x_2\}=\{x_2,x_3\}=0,\\
\{x_1,x_3\} &=-1,\ \{x_2,x_2\}=1
\endaligned\tag 4.4
$$
and also, 
$$
\re tz_0=x_1.\tag 4.5
$$
Observe that the conditions \thetag{4.4} are
equivalent to the fact that the (symmetric) matrix representation of
the bilinear form $\{\cdot,\cdot\}$ relative to the basis $x_1,x_2,x_3$
is by the matrix $\Lambda$, where
$$
\Lambda:=\pmatrix 0&0&-1\\0&1&0\\-1&0&0\endpmatrix.\tag 4.6
$$
We also write
$$\Lambda=(\lambda_{AB}).\tag 4.7$$
The set of all frames for
$\Gamma_\bC$ can be viewed as a real subgroup of the complex Lie
group $GL(\bC^4)$ as follows. If $(Z_0,Z_A)$ is a given frame,
then any other frame $(Z'_0,Z'_A)$ is obtained as
$$
(Z'_0,Z'_A)=(Z_0,Z_B)\pmatrix
v&0\\(k^B_1-v\delta^B_1)+iv^B&k^B_A\endpmatrix,\tag 4.8
$$
where $v,v^B\in \bR$ and the $3\times 3$-matrix $(k^B_A)$ satisfies
$$
k^C_Ak^D_B\lambda_{CD}=\lambda_{AB}\tag 4.9
$$
and also 
$$
\det(k^A_B)=1.\tag 4.10
$$
We denote by $H'$ the subgroup of $GL(\bC^4)$ consisting of all
matrices of the form
$$
\frak M=\pmatrix
v&0\\(k^B_1-v\delta^B_1)+iv^B&k^B_A\endpmatrix,\tag 4.11
$$
which satisfy
\thetag{4.9} and \thetag{4.10}. We also denote
by $K$ the subgroup
of $GL(\bR^3)$ which consists of
$(k^A_B)$ satisfying \thetag{4.9} and \thetag{4.10}. The group $K$
is isomorphic to the Lorenz group $SO(2,1)$. Indeed, if $O$ denote the
orthogonal
transformation for which $O^\tau \Lambda O$ equals the diagonal matrix
with $1,1,-1$ on the diagonal, then $K=O(SO(2,1))O^\tau$.

Note that \thetag{4.2}, \thetag{4.4}, and \thetag{4.5} imply that
$Z_0$, considered as the affine point $z_0\in\bC^3$ where $Z_0$ and
$z_0$ are as in \thetag{4.3}, can be viewed as
a point on $\Gamma_\bC$. We
denote by $H'_0$ the subgroup of $H'$ consisting of those matrices
which preserve $Z_0$ as an affine point on $\Gamma_\bC$, i.e. the
group of matrices of the form
$$
\frak  N=
\pmatrix
v&0\\0&k^B_A\endpmatrix,\tag 4.12
$$
where $(k^A_B)\in K$ with $k^A_1x_A=vx_1$. A straightforward
calculation shows that $(k^A_B)$ must be of the form
$$
(k^A_B)=\pmatrix v&a&\frac{1}{2}a^2v^{-1}\\0&1&av^{-1}\\
0&0&v^{-1}\endpmatrix,\tag 4.13
$$
for some $a\in\bR$. Thus, the group of frames for $\Gamma_\bC$ (which,
given a fixed frame, 
can be identified with the group $H'$ via \thetag{4.8}) is a principal
fiber bundle $P'\to \Gamma_\bC$ with group 
$H'_0$. Let us now choose an orientation, as explained in \S2, for the
Levi nullbundle of 
$\Gamma_\bC$ (which at a point $Z_0$ is spanned by $x_1$) and denote
by $P$ the group of frames consistent with this orientation. Then, as
is easily verified, $P$ is isomorphic to the group $H$, where
$H$ is the subgroup of $H'$ consisting of matrices of the form
\thetag{4.11} with $v>0$. If we also denote by $H_0$ the subgroup of
$H'_0$ consisting of those matrices of the form \thetag{4.12} for
which $v>0$, 
then $P\to M$ is a principal fiber bundle with group $H_0$. 
The reader should note that $\Gamma_\bC$ has two connected
components. This is reflected on the bundle $P$ by the fact that the
Lorenz group $SO(2,1)$ has two components.

A choice for the $4\times 4$-matrix $\Pi=(\pi^a_b)$ of
Maurer-Cartan forms for the group $H$ is given by
$$
dZ_a=\pi^b_a Z_b.\tag 4.14
$$
The Maurer-Cartan equations of structure then become
$$
d\pi^b_a=\pi^c_a\wedge\pi^b_c,\tag 4.15
$$
which follows directly from differentiating \thetag{4.14}.
Due to the form \thetag{4.3} of the frames $(Z_0,Z_A)$, we have
$\pi^0_0=dv/v$, $\pi^0_A=0$, and the $\pi^b_a$ are real. By
differentiating the defining equations \thetag{4.5} and using again
these equations, we deduce that the
$3\times 3$-matrix $(\pi^B_A)$ (which is a Maurer-Cartan matrix for
the group $K$)  is given by
$$
(\pi^B_A)=\pmatrix \pi^1_1&\pi^1_2&0\\
\pi^2_1&0&\pi^1_2\\
0&\pi^2_1&-\pi^1_1\endpmatrix,\tag 4.16
$$
for some real 1-forms $\pi^1_1,\pi^1_2,\pi^2_1$.
By differentiating $Z_0=v(1,z_0)$, we obtain
$$
dZ_0=\frac{dv}{v}(v,vz_0)+v(0,dz_0)=\frac{dv}{v}Z_0+(0,vdz_0).\tag 4.17
$$
Thus, by equation \thetag{4.5}, we have
$$
\frac{dv}{v}x_1+\frac{1}{2}(vdz_0+vd\bar z_0)=dx_1.\tag 4.18
$$
and hence, using also \thetag{4.14},
$$
\frac{1}{2}(\pi^A_0+\bar \pi^A_0)x_A=\left
(\pi^A_1-\frac{dv}{v}\delta^A_1\right) x_A.\tag 4.19
$$
Using \thetag{4.16}, we obtain the equations
$$
\left\{\aligned
&\pi^1_1-\pi^0_0=\frac{1}{2}(\pi^1_0+\bar \pi^1_0)\\
&\pi^2_1=\frac{1}{2}(\pi^2_0+\bar \pi^2_0)\\
&0=\pi^3_0+\bar \pi^3_0.\endaligned\right.\tag 4.20
$$
The last formula in \thetag{4.20} implies that $\pi^3_0$ is a purely
imaginary form. The matrix of 1-forms $\Pi=(\pi^a_b)$ is valued in the
Lie algebra $\frak h$ of $H$, and if we change frame using
\thetag{4.8} then the corresponding matrix of 1-forms $\Pi'$ for the
new frame $(Z'_0,Z'_A)$ is related to $\Pi$ by
$$
\Pi'=\ad(\frak M^{-1})\Pi:=\frak M^{-1}\Pi\frak M,\tag 4.21
$$
where $\frak M$ is the matrix given by \thetag{4.11}. It follows that
$\Pi$ is a Cartan connection on $P$ with group $H$ (see e.g.  [K, Chapter
IV] and also [CM]), which is flat (i.e. with vanishing curvature
form) by \thetag{4.15}.

Let us relate the above to the results in previous sections. It is
straightforward to verify that we can set
$$
\left\{\aligned
&\omega=-\frac{i}{2}\pi^3_0\\
&\omega^1=\frac{1}{2i}\pi^2_0\\
&\theta^2=\frac{1}{4i}\pi^1_0\\
&\Delta=\pi^0_0+\pi^1_1=2\pi^0_0+\frac{1}{2}(\pi^1_0+\bar\pi^1_0)\\
&\xi(=\Delta^1)=-\pi^1_2,\endaligned\right.\tag 4.22
$$
for the forms given by Theorem 3.1.37. Indeed, using \thetag{4.15} and
\thetag{4.20}, we obtain the equations
$$
\left\{\aligned
&d\omega=\Delta\wedge\omega+i\omega^{\bar 1}\wedge\omega^1\\
&d\omega^1=\xi\wedge\omega+\frac{1}{2}\Delta\wedge\omega^1-
i\omega^1\wedge\theta^2+i\theta^{\bar 2}\wedge\omega^1) +\omega^{\bar
1}\wedge \theta^2\\
&d\theta^2=\frac{1}{2}\xi\wedge\omega^1+2i\theta^{\bar 2}\wedge\theta^2\\
&d\Delta=i\xi\wedge\omega^1-i\xi\wedge\omega^{\bar 1}\\
&d\xi=-\frac{1}{2}(\Delta+2i\theta^2-2i\theta^{\bar
2})\wedge\xi,\endaligned \right.\tag 4.23
$$
which satisfy the conditions of Theorem 3.1.37. Note that the invariant
$\hat k$, defined in section 3, satisfies $\hat k\equiv 2i$. Thus, in
what follows, the phase function $t$ and the modulus $r$ as defined by
\thetag{3.1.2}, are identically 1. 

Recall that the group of
frames is a principal fiber bundle $P\to \Gamma_\bC$ with group
$H_0$. For any $\frak N\in H_0$,
where $\frak N\in H_0$ is given by \thetag{4.12} and $(k^A_B)$ by
\thetag{4.13}, a change of frame
$(Z'_0,Z'_A)=(Z_0,Z_B)\frak N$ results in the change of connection form
$$
\Pi'=\ad(\frak N^{-1})\Pi=\pmatrix
\pi^0_0&0\\v l^B_A\pi^A_0&l^B_A k^C_D\pi^A_C
\endpmatrix,\tag 4.24
$$
where $(l^C_D)$ denotes the inverse of $(k^A_B)$. Using
\thetag{4.22} and calculating the inverse of $(k^B_A)$ given by
\thetag{4.13}, we deduce that \thetag{4.24} yields the corresponding
transformation
$$\left\{\aligned
&\omega'=v^2\omega\\
&(\omega^1)'=-av\omega+v\omega^1\\
&(\theta^{2})'=\frac{1}{4}a^2\omega-\frac{1}{2}a\omega^1
+\theta^2\\  
&\Delta'=iav^{-1}(\omega^{\bar 1}-\omega^{1})+\Delta\\
&\xi'=\frac{1}{2}ia^2 v^{-1}(\omega^1-\omega^{\bar 1})-i
av^{-1}(\theta^2-\theta^{\bar
2})-\frac{1}{2}av^{-1}\Delta+v^{-1}\xi.\endaligned\right.\tag 4.25
$$
In particular, we obtain
$$
\pmatrix \omega'\\(\omega^1)'\\(\omega^{\bar 1})'\endpmatrix=\pmatrix
u&0&0\\ x&\sqrt{u}&0\\x&0&\sqrt{u}\endpmatrix\pmatrix
\omega\\{\omega^1}\\{\omega^{\bar 1}}\endpmatrix ,\tag 4.26
$$
where
$$
u=v^{2},\quad x=-av.\tag 4.27
$$
Hence, we have defined an isomorphism $\phi\:H_0\to G_2$, where $G_2$
is as defined in
section 3.1 with $t\equiv 1$, defined by
$$
\phi\pmatrix v&0&0&0\\0&v&a&\frac{1}{2}a^2v^{-1}\\0&0&1&av^{-1}\\0&
0&0&v^{-1}\endpmatrix=\pmatrix v^{2}&0&0\\
-av&v&0\\-av&0&v\endpmatrix .\tag 4.28
$$

Let $M$ be any $5$-dimensional CR manifold of
hypersurface type which satisfies the conditions in Theorem
3.1.37. Furthermore, we assume that the invariant $\hat k\equiv 2i$ in
a neighborhood of $p_0$. Let $\omega,\omega^1,\theta^2,\Delta,\xi$ be
the forms given 
by Theorem 3.1.37 such that $\underline{\omega}$ is given by
\thetag{3.1.36}. Define the $\frak h$-valued 1-form $\Pi=(\pi_a^b)$
by 
$$
\Pi:=\pmatrix \pi^0_0&0&0&0\\\pi^1_0&\pi^1_1&\pi^1_2&0\\\pi^2_0&
\pi^2_1&0&\pi^1_2\\\pi^3_0&
0&\pi^2_1&-\pi^1_1\endpmatrix,\tag 4.29
$$
where $$\pi_0^0,\pi^1_0,\pi^2_0,\pi^3_0,\pi^1_1,\pi^2_1,\pi^1_2,$$ are
obtained 
by solving \thetag{4.22}, using also the first two equation of
\thetag{4.20}. 
Clearly, $\Pi$ defines an isomorphism between $T_y Y_2$, where $Y_2\to
M$ is the principal bundle given by Theorem 3.1.37, and $\frak
h$ for every $y\in Y_2$. However, it is not difficult to verify that
$\Pi$ is in general 
not a Cartan connection, i.e.\ it does not transform according to
\thetag{4.21}. Nevertheless, by defining the curvature
$$
\Omega=d\Pi-\Pi\wedge\Pi,\tag 4.30
$$ 
a direct consequence of Cartan's solution of the equivalence problem
for $\{1\}$- structures (see e.g.\ [G]) is the following
characterization of the tube over the light cone.

\proclaim{Theorem 4.31} Let $M$ be a $5$-dimensional real-analytic CR
manifold of hypersurface type which is $2$-nondegenerate and Levi
uniform of rank $1$ at $p_0\in M$. Assume that $\hat k\equiv 2i$. 
Choose an orientation for the Levi 
nullspace $\frak N_{p_0}$ (as explained in \S$2$), and denote by
$Y_2\to M$ the principal 
bundle with $1$-form $\underline{\omega}$ given by Theorem
$3.1.37$. Then, the $\frak 
h$-valued $1$-form $\Pi$ defined by \thetag{4.29}, where
$$
\left\{\aligned
&\pi^0_0:=\frac{1}{2}(\Delta-2i\theta^2+2i\theta^{\bar 2})\\
&\pi^1_1:=\frac{1}{2}(\Delta+2i\theta^2-2i\theta^{\bar 2})\\
&\pi^1_0:=4i\theta^2\\
&\pi^2_0:=2i\omega^1\\
&\pi^3_0:=2i\omega\\
&\pi^2_1:=i(\omega^1-\omega^{\bar 1})\\
&\pi^1_2:=-\xi\endaligned\right.\tag 4.32
$$
and $\underline{\omega}$ is given by \thetag{3.1.36}, defines an 
isomorphism $T_y Y_2\cong 
\frak h$ for every $y\in Y_2$, with the following property. 
There exists a local real-analytic CR diffeomorphism $f\:M\to
\Gamma_\bC$ near $p_0$ if and only if the curvature $\Omega$ given by
\thetag{4.30} vanishes identically. 
\endproclaim

We conclude the discussion of the tube over the light cone by 
computing the dimension of the stability group $\Aut(\Gamma_\bC,p_0)$ of 
$\Gamma_\bC$ at a 
point $p_0\in \Gamma_\bC$. Observe that, given a frame $(Z_0,Z_A)$ in
$P$, the manifold $\Gamma_\bC$ can be viewed as the
quotient group $H/H_0$ via the identification $P\cong H$ provided by
\thetag{4.8}. Let us denote the affine point on $\Gamma_\bC$
corresponding to $Z_0$ by $p_0$. Then, under the identification
$\Gamma_\bC\cong H/H_0$, $p_0$ corresponds to the coset
$eH_0$, where $e\in H$ denotes the identity matrix. The group $H_0$
acts on the left on $H/H_0$ and each homomorphism $aH_0\mapsto baH_0$,
for $b\in H_0$, 
preserves the point $p_0\cong eH_0$. It is straightforward to verify that
the action is effective; i.e. if, for $b\in H_0$, the homomorphism
$aH_0\mapsto baH_0$ is the identity, then $b=e$. Let us denote by
$f_b\:(\Gamma_\bC, p_0)\to (\Gamma_\bC,p_0)$ the mapping corresponding
to the homomorphism $aH_0\mapsto baH_0$. Each $f_b$ is a CR
diffeomorphism. (Indeed, it is not difficult to compute $f_b$ in
coordinates and see that $f_b$ is induced by an invertible linear
transformation of $\bC^3$.) Thus, $b\mapsto f_b$ embeds $H_0$ as a subgroup
of $\Aut(\Gamma_\bC,p_0)$. Since $\dim H_0=2$, we conclude that $\dim
\Aut(\Gamma_\bC,p_0)\geq 
2$. On the other hand, by Theorem 3.1.37 and 
[K, Theorem 3.2], it follows (as in the introduction) that the
subgroup of $\Aut(\Gamma_\bC, p_0)$ consisting of those CR
diffeomorphisms that preserve the orientation of the Levi nullspace
chosen above embeds as a closed submanifold of $P_{p_0}\cong
H_0$. Hence, we have $\dim\Aut(\Gamma_\bC,p_0)=2$.

%
%

\heading 5. Concluding remarks; the higher dimensional
case\endheading

Let us briefly return to the situation in section 2, i.e.\ $M$ is a
smooth CR manifold (of hypersurface type and dimension $2n+1$) which
is Levi uniform of rank $n-1$ at $p_0$. We also assume that $M$
satisfies Condition 2.21 and 2.25 (which in particular imply that $M$
is pseudoconvex and $2$-nondegenerate at $p_0$). As in section 3, we
consider 
equation \thetag{2.35} which, in view of \thetag{2.13}, can be
rewritten as follows 
$$
\tilde k_{\bar\mu}=\hat k_{\bar \mu}+2i\hat g^{\bar\alpha\nu}
c_\nu\hat h_{\bar\alpha\bar\mu}-\overline{c_{\mu}}\hat k.\tag 5.1
$$
Now, using \thetag{2.29} we deduce that
$$
\tilde k_{\bar\mu}=\hat k_{\bar \mu}+2i
c_\mu\lambda_{\bar\mu}-\overline{c_{\mu}}\hat k.\tag 5.2
$$
Thus, either $|\hat k(p_0)|\neq 2\lambda_{\bar\mu}$ for 
$\mu=1,2,\ldots, n-1$, or $|\hat k(p_0)|\neq 2\lambda_{\bar\mu_0}$ for
some $\mu_0\in\{1,2,\ldots, n-1\}$. In the first case, we can solve
for each $c_\mu$ in the equation $\tilde k_{\bar\mu}=0$ and proceed as
in section 3.2 to construct a principal bundle $P\to M$ with 1-form
$\underline{\omega}$ reducing the CR structure on $M$ to a
parallelism. In the latter case, we can solve for $c_{\bar\mu}$ in
the equation $\tilde k_{\bar \mu}=0$ for all $\mu\neq\mu_0$. We then
proceed as in in section 3.1 to determine $c_{\bar\mu_0}$ and
construct the bundle $P\to M$ with 1-form $\underline{\omega}$. We do
not give the details here. Conditions 2.21 and
2.25 do not appear to be natural when $n\geq 3$. (Recall,
however, that 
these two conditions reduce to $2$-nondegeneracy and Levi uniformity
when $n=2$.) In
particular, the tube over the light cone in $\bC^{n+1}$, $n\geq 3$,
does not satisfy these conditions.

\enddocument
\end